\documentclass{amsart}

\usepackage{amssymb, amsbsy, amsthm, amsmath, amstext, amsopn, verbatim}
\usepackage[all]{xy}
\usepackage{amsfonts}
\usepackage{amscd}
\newtheorem{thm}{Theorem} [section]
\newtheorem{lemma}[thm]{Lemma}

\newtheorem{corollary}[thm]{Corollary}
\newtheorem{prop}[thm]{Proposition}

\theoremstyle{definition}

\newtheorem{defn}[thm]{Definition}

\newtheorem{construction}[thm]{Construction}

\theoremstyle{remark}

\newtheorem{remark}[thm]{Remark}

\begin{document}

\numberwithin{equation}{section}

\newcommand{\hs}{\mbox{\hspace{.4em}}}
\newcommand{\ds}{\displaystyle}
\newcommand{\bd}{\begin{displaymath}}
\newcommand{\ed}{\end{displaymath}}
\newcommand{\bcd}{\begin{CD}}
\newcommand{\ecd}{\end{CD}}

\newcommand{\on}{\operatorname}
\newcommand{\proj}{\operatorname{Proj}}
\newcommand{\Proj}{\operatorname{Proj}}
\newcommand{\bproj}{\underline{\operatorname{Proj}}}
\newcommand{\spec}{\operatorname{Spec}}
\newcommand{\bspec}{\underline{\operatorname{Spec}}}
\newcommand{\pline}{{\mathbf P} ^1}
\newcommand{\aline}{{\mathbf A} ^1}
\newcommand{\pplane}{{\mathbf P}^2}
\newcommand{\aplane}{{\mathbf A}^2}
\newcommand{\coker}{{\operatorname{coker}}}
\newcommand{\ldb}{[[}
\newcommand{\rdb}{]]}

\newcommand{\Sym}{\operatorname{Sym}^{\bullet}}
\newcommand{\Symp}{\operatorname{Sym}}
\newcommand{\Pic}{\operatorname{Pic}}
\newcommand{\AAut}{\operatorname{Aut}}
\newcommand{\PAut}{\operatorname{PAut}}

\newcommand{\too}{\twoheadrightarrow}
\newcommand{\Z}{{\mathbf Z}}
\newcommand{\C}{{\mathbf C}}
\newcommand{\Cx}{{\mathbf C}^{\times}}
\newcommand{\cA}{{\mathcal A}}
\newcommand{\cS}{{\mathcal S}}
\newcommand{\cV}{{\mathcal V}}
\newcommand{\cM}{{\mathcal M}}
\newcommand{\bA}{{\mathbf A}}
\newcommand{\cB}{{\mathcal B}}
\newcommand{\cC}{{\mathcal C}}
\newcommand{\cD}{{\mathcal D}}
\newcommand{\D}{{\mathcal D}}
\newcommand{\cs}{{\mathbf C} ^*}
\newcommand{\boldc}{{\mathbf C}}
\newcommand{\cE}{{\mathcal E}}
\newcommand{\cF}{{\mathcal F}}
\newcommand{\bF}{{\mathbf F}}
\newcommand{\cG}{{\mathcal G}}
\newcommand{\G}{{\mathbf G}}
\newcommand{\cH}{{\mathcal H}}
\newcommand{\cJ}{{\mathcal J}}
\newcommand{\cK}{{\mathcal K}}
\newcommand{\cL}{{\mathcal L}}
\newcommand{\baL}{{\overline{\mathcal L}}}
\newcommand{\M}{{\mathcal M}}
\newcommand{\bM}{{\mathbf M}}
\newcommand{\bm}{{\mathbf m}}
\newcommand{\cN}{{\mathcal N}}
\newcommand{\theo}{\mathcal{O}}
\newcommand{\cP}{{\mathcal P}}
\newcommand{\cR}{{\mathcal R}}
\newcommand{\boldp}{{\mathbf P}}
\newcommand{\boldq}{{\mathbf Q}}
\newcommand{\bbL}{{\mathbf L}}
\newcommand{\cQ}{{\mathcal Q}}
\newcommand{\cO}{{\mathcal O}}
\newcommand{\Oo}{{\mathcal O}}
\newcommand{\OX}{{\Oo_X}}
\newcommand{\OY}{{\Oo_Y}}
\newcommand{\otY}{{\underset{\OY}{\ot}}}
\newcommand{\otX}{{\underset{\OX}{\ot}}}
\newcommand{\cT}{{\mathcal T}}
\newcommand{\cU}{{\mathcal U}}\newcommand{\cX}{{\mathcal X}}
\newcommand{\cW}{{\mathcal W}}
\newcommand{\boldz}{{\mathbf Z}}
\newcommand{\qgr}{\operatorname{q-gr}}
\newcommand{\gr}{\operatorname{gr}}
\newcommand{\rk}{\operatorname{rk}}
\newcommand{\coh}{\operatorname{coh}}
\newcommand{\End}{\operatorname{End}}
\newcommand{\Hom}{\operatorname{Hom}}
\newcommand{\uHom}{\underline{\operatorname{Hom}}}
\newcommand{\uHomY}{\uHom_{\OY}}
\newcommand{\uHomX}{\uHom_{\OX}}
\newcommand{\Ext}{\operatorname{Ext}}
\newcommand{\bExt}{\operatorname{\bf{Ext}}}
\newcommand{\Tor}{\operatorname{Tor}}

\newcommand{\inv}{^{-1}}
\newcommand{\airtilde}{\widetilde{\hspace{.5em}}}
\newcommand{\airhat}{\widehat{\hspace{.5em}}}
\newcommand{\nt}{^{\circ}}
\newcommand{\del}{\partial}

\newcommand{\supp}{\operatorname{supp}}
\newcommand{\GK}{\operatorname{GK-dim}}
\newcommand{\hd}{\operatorname{hd}}
\newcommand{\id}{\operatorname{id}}
\newcommand{\res}{\operatorname{res}}
\newcommand{\lrar}{\leadsto}
\newcommand{\im}{\operatorname{Im}}
\newcommand{\HH}{\operatorname{H}}
\newcommand{\TF}{\operatorname{TF}}
\newcommand{\Bun}{\operatorname{Bun}}
\newcommand{\BunD}{\operatorname{Bun}_{\D}}

\newcommand{\PicD}{\operatorname{Pic}_{\D}}
\newcommand{\Hilb}{\operatorname{Hilb}}
\newcommand{\Fact}{\operatorname{Fact}}
\newcommand{\CM}{\mathfrak{CM}}
\newcommand{\ECM}{\mathfrak{ECM}}
\newcommand{\HCM}{\mathfrak{HCM}}
\newcommand{\MD}{\mathfrak{M}^{\D}}
\newcommand{\F}{\mathcal{F}}
\newcommand{\Ff}{\mathbb{F}}
\newcommand{\nthord}{^{(n)}}
\newcommand{\Aut}{\underline{\operatorname{Aut}}}
\newcommand{\Gr}{{\on{Gr}}}
\newcommand{\GR}{\operatorname{GR}}
\newcommand{\GRo}{{\operatorname{GR}^{\circ}}}
\newcommand{\GRon}{\operatorname{GR}^{\circ}_n}
\newcommand{\Fr}{\operatorname{Fr}}
\newcommand{\GL}{\operatorname{GL}}
\newcommand{\gl}{\mathfrak{gl}}
\newcommand{\SL}{\operatorname{SL}}
\newcommand{\ff}{\footnote}
\newcommand{\ot}{\otimes}
\def\Ext{\operatorname {Ext}}
\def\Hom{\operatorname {Hom}}
\def\Ind{\operatorname {Ind}}
\def\bbZ{{\mathbb Z}}
\newcommand{\mOp}{\on{MOp}}

\newcommand{\AutO}{\on{Aut}\Oo}
\newcommand{\Der}{{\on{Der}\,}}
\newcommand{\DerO}{{\on{Der}\Oo}}
\newcommand{\AutK}{\on{Aut}\K}
\newcommand{\SShv}{\on{SShv}}

\newcommand{\nc}{\newcommand}
\nc{\Mf}{{\mathfrak M}} \nc{\cont}{\on{cont}} \nc{\rmod}{\on{mod}}
\nc{\Mtil}{\widetilde{M}} \nc{\ol}{\overline} \nc{\wb}{\overline}
\nc{\wt}{\widetilde} \nc{\wh}{\widehat} \nc{\sm}{\setminus}
\nc{\mc}{\mathcal} \nc{\mbb}{\mathbb} \nc{\Mbar}{\wb{M}}
\nc{\Nbar}{\wb{N}} \nc{\Mhat}{\wh{M}} \nc{\pihat}{\wh{\pi}}
\nc{\JYX}{\cJ_{Y\leftarrow X}} \nc{\phitil}{\wt{\phi}}
\nc{\Qbar}{\wb{Q}} \nc{\DYX}{\D_{Y\leftarrow X}} \nc{\DXY}{\D_{X\to
Y}} \nc{\dR}{\stackrel{\bbL}{\underset{\D_X}{\ot}}}
\nc{\Winfi}{\cW_{1+\infty}} \nc{\K}{{\mc K}} \nc{\Kx}{{\mc
K}^{\times}} \nc{\Ox}{{\mc O}^{\times}} \nc{\unit}{{\bf \on{unit}}}
\nc{\boxt}{\boxtimes} \nc{\xarr}{\stackrel{\rightarrow}{x}}
\nc{\Gamx}{\Gamma_{-}^{\times}} \nc{\Gap}{\Gamma_+}
\nc{\Gamxn}{\Gamma_{-,n}^{\times}} \nc{\Gapn}{\Gamma_{+,n}}

\nc{\Dlog}{{\mathcal D}_{\on{log}}} \nc{\Ereg}{E_{sm}}
\nc{\Jac}{\on{Jac}} \nc{\Vplus}{V^{+}} \nc{\Vc}{V^{\vee}}
\nc{\af}{{\mathfrak a}} \nc{\h}{{\mathfrak h}}
\nc{\BunDP}{\on{Bun}_{\boldp(\cD)}} \nc{\SSh}{\on{SShv}}
\nc{\Enat}{E^{\natural}} \nc{\Enatbar}{\overline{E}^{\natural}}
\nc{\Id}{\on{Id}} \nc{\Fin}{F^{-1}} \nc{\Dhat}{\wh{D}}
\nc{\Of}{{\mathbb{O}}}\nc{\Pp}{{\mbb P}} \nc{\Hitch}{{\on{Hitchin}}}

\nc{\framing}{T} \nc{\At}{A} \nc{\ul}{\underline}
\nc{\Dx}{D^{\times}}

\title{Flows of Calogero-Moser Systems}
\author{David Ben-Zvi}
\address{Department of Mathematics\\University of Texas at Austin\\Austin, TX 78712 USA}
\email{benzvi@math.utexas.edu}
\author{Thomas Nevins}
\address{Department of Mathematics\\University of Illinois at Urbana-Champaign\\Urbana, IL 61801 USA}
\email{nevins@uiuc.edu}


\maketitle

\section{Introduction}
The Calogero-Moser (or CM)
particle system \cite{Ca1, Ca2} 
and its generalizations appear, in a variety of ways, in
 integrable systems, nonlinear PDE, representation theory,
and string theory.  Moreover, the partially completed
CM systems---in which
dynamics of particles are continued through collisions---have been
identified as meromorphic Hitchin systems (see, for example, \cite{BBT, DM, GorNe,HuMark, HN, Kr2, Kr Lax Vect, Nekrasov} among others).
This Hitchin-type description gives natural
``geometric action-angle variables'' for the CM system.

Motivated by relations of the CM system to nonlinear PDE
\cite{announce, solitons}, we  introduce a new class of
generalizations of the spin CM particle systems, the {\em framed}
(rational, trigonometric and elliptic) CM systems. We give two
algebro-geometric descriptions of these systems: the first, perhaps more
familiar, description
uses meromorphic $GL_n$ Hitchin systems with decorations
(framing data) on (cuspidal, nodal and smooth) cubic curves.  The second
description identifies these phase spaces with moduli spaces of one-dimensional
sheaves on corresponding ``twisted'' ruled surfaces.

We also present a simple
geometric formulation of the flows of all meromorphic $GL_n$ Hitchin
systems (with no regularity assumptions) as {\em tweaking} flows on
spectral sheaves. Using this formulation, we show that all spin and framed CM
systems are identified with hierarchies of tweaking flows on the
corresponding spectral sheaves. This generalizes the well-known
description of spinless CM systems in terms of tangential covers
(see \cite{TV,TV2}). In \cite{solitons}, we prove that a Fourier transform
identifies
 the framed CM systems (in their
spectral incarnation) as the particle systems describing
the motion of poles of all meromorphic (rational, trigonometric, and
elliptic) solutions of (generalized) multicomponent KP hierarchies.

We begin with an overview of the paper. Sections \ref{CM section}
and \ref{CM and Hitchin chapter} review relevant background, Section
\ref{CM spectral sheaves} introduces framed CM systems and the
corresponding spectral sheaves on ruled surfaces, and Section
\ref{flows section} contains the description of flows of Hitchin and
framed CM systems.

\subsection{Review of Calogero-Moser Systems}
In Section \ref{CM section}, we review the definition of the
Calogero-Moser (CM) particle systems. The Calogero-Moser systems
are a family of completely integrable classical hamiltonian systems,
describing particles on the line interacting with a quadratic
potential. The systems come in three basic variants: rational,
trigonometric and elliptic, according to the type of function used
in the potential. We will work with complexified Calogero-Moser
systems, in which the three variants naturally correspond to systems
of interacting particles on a one-dimensional complex group $\G$,
namely $\C$, $\Cx$ or an elliptic curve. The Calogero-Moser systems
also have natural partial completions, in which we allow the
positions of particles to collide \cite{KKS, Wilson CM}. The
(completed) Calogero-Moser systems feature in a variety of
seemingly unrelated areas, in particular as describing the motion of
poles of meromorphic solutions to KdV and KP equations (see
\cite{AMM,CC,Kr1,Kr2,TV2,Wilson CM} and \cite{Benn,BAMS,Wsurv3} for
surveys), and as analogs of Hilbert schemes in noncommutative
geometry \cite{BW ideals,BGK1,KKO,BraNe}.

The Calogero-Moser systems have a natural generalization, in which
the particles carry spins, which are covectors and vectors $u_i,v_i$
in an auxiliary vector space $\C^k$ \cite{GH}. The dynamics of the
spin CM systems depends on the spins only through the pairings
$f_{ij}=u_i(v_j)$, so one usually considers instead a reduced
version (the {\em Euler-Calogero-Moser system}) in which we only
keep track of the $f_{ij}$ \cite{BBKT,Res,Nekrasov}. From the point
of view of noncommutative geometry and KP equations (e.g.
\cite{BGK2,solitons}), it is the full phase space of the spin CM
hierarchy, and its framed generalizations presented in this paper,
which play a central role. We review the definition of the completed
rational and trigonometric spin CM systems in terms of quiver data
(pairs of matrices, vectors and covectors satisfying a shifted
moment map condition) in Section \ref{CM matrices}.

\subsection{Calogero-Moser and Hitchin Systems}
We begin Section \ref{CM and Hitchin chapter} by reviewing (in
Section \ref{cubics}) the basic setup, the family of Weierstrass
cubic curves $E$ (irreducible genus one curves). These come in three
variants, according to the type of the smooth locus $\G\subset E$:
rational (the cuspidal cubic, with $\G=\C$), trigonometric (the
nodal cubic, with $\G=\Cx$), and elliptic ($\G=E$ is an elliptic
curve).

In Section \ref{CM and Hitchin} we review the relation between
meromorphic Hitchin systems\cite{Hi,M,DM} on cubic curves and
Euler-Calogero-Moser systems (following \cite{Nekrasov, GorNe} as
well as the reviews \cite{DM,BBT}). The Hitchin formulation is a
variant of the Lax form with spectral parameter for elliptic CM
systems, due to Krichever \cite{Kr2} (see also \cite{Kr Lax Vect}
where elliptic CM systems appear as part of a general hamiltonian
theory of Lax operators on algebraic curves).

Recall that the ($GL_n$) Hitchin systems are integrable systems on
moduli spaces of vector bundles on a curve, equipped with a Higgs
field (endomorphism-valued one-form). A key ingredient in this
approach to particle systems is the identification
(\cite{FriedMorgminuscule}, see Section \ref{particles as bundles})
between positions of particles (configurations of points, or more
generally torsion sheaves, on $\G$) and vector bundles on $E$, which
is a special case of the Fourier-Mukai transform on cubic curves
\cite{BuK,solitons}.

\subsection{Framed Calogero-Moser Systems and Spectral Sheaves}
Section \ref{CM spectral sheaves} introduces and studies the framed
Calogero-Moser systems. These systems are best described
geometrically using a nontrivial affine bundle (the Serre surface)
$\Enat\to E$ modelled on $T^*E$. This formulation for the elliptic
(spinless) Calogero-Moser system was explored in detail by
Treibich-Verdier in their work on tangential covers \cite{TV,TV2},
and in \cite{DW,Donagi} (see \cite{DM}). One of our goals is to
generalize this description to include the rational and
trigonometric cases and to CM particles with spins or {\em
framings}.

The ruled surface $\Enatbar\to E$ over a cubic curve $E$ is
discussed in Section \ref{about Enat}. It comes equipped with a
section $E_{\infty}$, whose complementary affine bundle
$\Enat=\Enatbar\setminus E_{\infty}$ is the natural home for the
Weierstrass $\zeta$-function. The surface $\Enat$ is naturally
birational to $T^*E$ (again using the $\zeta$-function), so that the
complement of the fibers over the basepoint $b\in E$ (identity
element of the group $\G$) are identified.

In Section \ref{twisted Higgs} we introduce the notion of twisted
Higgs field, which relate to sheaves on $\Enat$ in the same way that
ordinary Higgs fields relate to sheaves on $T^*E$. Using the
birational identification of $\Enat$ and $T^*E$ we have a simple
bijection between regular and twisted (framed or unframed) Higgs
fields, which shifts the Higgs field by $\zeta\on{Id}$. This shift
provides a geometric origin (the transition from $T^*E$ to $\Enat$)
for the appearance of the shift by $\on{Id}$ in the hamiltonian
reduction and Hitchin system descriptions of CM hamiltonians (and
corresponds, under the extended Fourier-Mukai transform
\cite{solitons}, to a transition from sheaves on $T^*E$ to
$\D$-modules).

For any torsion coherent sheaf $\framing$ on the group $\G\subset
E$, we introduce in Section \ref{framed CM systems} the
$\framing$-{\em framed Calogero-Moser systems}, which are Hitchin
systems on $\framing$-{\em framed} Higgs (or twisted Higgs) bundles.
The framing consists of a factorization of the polar parts of a
Higgs field into a map to $\framing$. This generalizes the usual
$k$-spin Calogero-Moser systems, which are recovered when the
framing $\framing=\C^k$, considered as a sum of skyscrapers
$\Oo_b^{\oplus k}$ at the basepoint of $\G$. The framing is also
analogous to that appearing in the definition of Nakajima quiver
varieties \cite{Nakajima}. Thanks to the identification between
vector bundles on $E$ and configurations of particles, these systems
have the form of generalized spin particle systems.

In Section \ref{CM spectral sheaves section} we introduce the spaces
of {\em framed CM spectral sheaves}, which provide a geometric phase
space for spin (and framed) CM systems. Recall that the phase spaces
of Hitchin systems on a curve $X$ are described geometrically as
spaces of spectral sheaves on the cotangent bundle $T^*X$---i.e.
line bundles (or more generally torsion-free sheaves) supported on
curves in the surface $T^*X$.
 In fact the description of an
integrable system by spectral sheaves is a geometric version of
transforming the system in action-angle variables---the support
curve is invariant under the system and plays the role of the action
variables, while the line bundle on the curve plays the role of the
angle variable.

For CM systems, the natural spectral sheaves live on the ruled
surface $\Enatbar$. Given a torsion sheaf $\framing$ on the smooth
locus $\G\subset E$, we define a {\em $\framing$-framed CM spectral
sheaf} on $E$ to be a torsion-free sheaf supported on a curve in
$\Enatbar$, whose restriction to $E_\infty$ is identified with
$\framing$---see Definition \ref{CM spectral sheaves def} for a
precise definition (the algebraic geometry of related linear series
on $\Enatbar$ is studied in \cite{T matrix}). The torsion sheaf
$\framing$ again plays the role of the spin variables of the
corresponding Calogero-Moser particles---in particular we'll show
(Corollary \ref{CM is CM}) that the $k$-spin CM system is realized
as the case where $\framing$ is the vector space $\C^k$, considered
as a skyscraper sheaf $\Oo_b^{\oplus k}$ at the basepoint $b\in E$.

In Section \ref{spectral sheaves and higgs} we identify moduli
spaces of $\framing$-framed CM spectral sheaves and Higgs bundles:

\begin{thm}[Theorem \ref{CM identification}] There is a canonical
isomorphism $\CM_n(E,\framing) \to \HCM_n(E,\framing)$ between the
moduli of $\framing$-framed spectral sheaves and $\framing$--framed
Higgs fields.
\end{thm}

The identification is based on a result of Katzarkov, Orlov and
Pantev \cite{KOP} which uses Koszul duality to identify moduli
spaces of framed sheaves on ruled surfaces over curves in terms of
linear algebra data on the curve. It turns out that the Koszul data
for CM spectral sheaves are precisely framed twisted Higgs bundles,
so that the Higgs field gives the structure of the underlying
spectral sheaf while the factorization of the poles into spin
variables gives the framing of the spectral sheaf. In particular,
when $\framing=\Oo_b^k$, we find that the phase spaces of the
rational, trigonometric and elliptic spin CM particles are
identified with spaces of framed spectral sheaves on $\Enat$.

\subsection{Flows on Spectral Sheaves}
In Section \ref{flows section}, the heart of the paper, we study
flows on moduli spaces of spectral sheaves. We first discuss in
Section \ref{tweaking section} the general principle (tweaking),
whereby sheaves are deformed by an infinitesimal version of
tensoring with line bundles. Explicitly we can construct
deformations of arbitrary sheaves from germs of meromorphic
functions. Namely, multiplication by such a function gives a germ of
a meromorphic endomorphism of any sheaf, using which we ``change the
transition functions" (or define an $\Ext^1$ class).

In Section \ref{Hitchin flows section} we present a geometric
description of the flows of meromorphic $GL_n$ Hitchin systems on an
algebraic curve $X$ as tweaking flows (Theorem \ref{Hitchin post}).
(Related geometric pictures with various restrictions on the allowed
spectral curves appear for example in \cite{AB,AHP,DM,BF,LiMulasePrym, LiMulaseHitchin}.) Namely we
point out an obvious bijection between Hitchin hamiltonians for
vector bundles on a curve $X$ and classes in $\HH^1(T^*X,\Oo)$, and
likewise for meromorphic Hitchin systems and meromorphic germs on
$T^*X$. It is then an easy check that the corresponding tweaking and
hamiltonian flows agree (this is a global analog of the trivial
spectral description of the hamiltonian flows of trace polynomials
on $T^*GL_n$).  We thus realize the Hitchin flows not in terms of
the action of line bundles on each specific spectral curve but
uniformly as the infinitesimal version of the action of line bundles
on $T^*X$ (i.e. as the action of a commutative Lie {\em algebra},
rather than Lie algebroid).

\begin{thm}[Theorem \ref{Hitchin post}]
The hamiltonian flows on the moduli of (possibly meromorphic) Higgs
bundles on a curve $X$ are given by the multiplication action of
classes in $\HH^1(T^*X,\Oo)$ (or generally meromorphic germs on
$T^*X$) on the corresponding spectral sheaves.
\end{thm}

This concrete realization of Hitchin flows has various applications.
As an example we describe, as Corollary \ref{compatibility}, a
simple generalization of the Compatibility Theorem of Donagi-Markman \cite{DM}
and Li-Mulase \cite{LiMulaseHitchin}
relating Hitchin and Heisenberg (or KdV) flows,
dropping all assumptions on the regularity of the Higgs field.

Finally, in Section \ref{CM flows} we define the framed
Calogero-Moser hierarchies, as tweaking flows on framed spectral
sheaves. Specifically, we tweak CM spectral sheaves by principal
parts of functions with poles along the curve $E_\infty$ at
infinity. We are then able to identify explicitly all of these flows
with Hitchin hamiltonian flows, and in particular with the explicit
form of the spin (and framed) Calogero-Moser hamiltonians on
particles. We summarize as follows:

\begin{thm}[Theorem \ref{framed CM flows
thm}]
\mbox{}
\begin{enumerate}
\item The flows of the rational, trigonometric and elliptic spin CM hamiltonian
systems are identified with explicit tweaking flows (see Definition
\ref{framed CM hierarchy def}) along $E_\infty$
on ($\Oo_b^k$--framed) CM spectral sheaves on $\Enat$.
\item More generally, for 
any torsion sheaf $\framing$ on $\G\subset E$ with $E$ a cubic curve,
the flows of the $\framing$-framed Calogero-Moser hamiltonians are
identified (under the isomorphism
$\CM_n(E,\framing)\to\HCM_n(E,\framing)$ above) with explicit
tweaking flows (Definition
\ref{framed CM hierarchy def}) along $E_\infty$ on $\framing$-framed CM spectral
sheaves.
\item For simple framing $\framing$ (for example in the spinless case $\framing=\Oo_b$),
the CM hamiltonian flows generate all tweaking flows along $E_\infty$.
\end{enumerate}
\end{thm}

\subsection{Motivation} The impetus for the present work comes from
the correspondence between meromorphic solutions of soliton
equations of KP type and particle systems of CM type. In
\cite{solitons} (see \cite{announce} for an overview) we establish a
very broad form of this correspondence, generalizing and refining (in
particular) the results of \cite{Kr2,Wilson CM,TV} in the spinless
case. Namely, we apply an extension of the Fourier-Mukai transform
to the spaces of framed CM spectral sheaves studied in this paper,
extending the identification of the underlying vector bundles with
configurations of particles. We obtain an identification of these
phase spaces with moduli spaces of {\em framed $\D$-bundles} on $E$.
These moduli spaces are noncommutative analogs of Hilbert schemes
(in the rank one case) or framed torsion--free sheaves on $T^*E$,
and the isomorphism may be considered a separation of variables \`a
la Sklyanin (see \cite{GNR}) for elliptic Hitchin systems. This
generalizes the relation between rational Calogero-Moser spaces and
ideals in the Weyl algebra \cite{LB, BW ideals, BGK1}. Moreover, we
show that framed $\D$-bundles provide a natural geometric phase
space for the meromorphic (rational, trigonometric and elliptic)
multicomponent KP hierarchy, and that the isomorphism of moduli
spaces identifies the KP and CM flows. The positions of the CM
particles are identified with the ``singularities" of the
$\D$-bundles, which are the poles of the corresponding meromorphic
KP solution. Thus, framed CM systems describe the motion of poles of
general meromorphic solutions of multicomponent KP hierarchies.

\subsection{Acknowledgments}
The authors are grateful to R. Donagi, N. Nekrasov and T. Pantev for
helpful conversations. In particular T. Pantev explained to us his
work with Katzarkov and Orlov describing framed sheaves by Koszul
data which we use to identify Calogero-Moser spectral sheaves with
Higgs bundles. The first author would also like to express his
appreciation to V. Drinfeld: it was while speaking in Drinfeld's seminar
that the first author
first learned about Calogero-Moser systems and the role of cubic
curves in representation theory.

The first author was supported in part by an
 NSF postdoctoral fellowship at the University of Chicago and an
MSRI postdoctoral fellowship, as well as by NSF CAREER grant
DMS-0449830 at the University of Texas. The second author was
supported in part by an NSF postdoctoral fellowship at the
University of Michigan and an MSRI postdoctoral fellowship, as well
as by NSF grant DMS-0500221.

\section{Review of Calogero-Moser Systems}\label{CM
section}

\subsection{Introducing the Spin Calogero-Moser System}

In this section we discuss the spin Calogero-Moser system,
following \cite{GH,BBKT,Nekrasov,Res}---see also the chapter in
\cite{BBT}. See \cite{announce} for a review of the usual (spinless)
complexified Calogero-Moser system following \cite{KKS,Wilson
CM,Nekrasov,GorNe,Enriquez,Martinec}.

Connected one-dimensional complex algebraic groups $\G$ fall into
three classes: the additive group $\C$, the multiplicative group
$\Cx$, and the one-parameter family of elliptic curves $E$. These
cases fall under the monikers {\em rational, trigonometric} and {\em
elliptic} according to the type of functions on the universal cover
$\C$ which correspond to meromorphic functions on $\G$.

The $k$-spin $n$-particle Calogero-Moser system is a hamiltonian
system describing $n$ identical particles on $\G$ equipped with
spins in the auxiliary $k$-dimensional vector space $\C^k$. Thus,
consider $n$ distinct points (positions) $q_1, \dots, q_n$ in $\G$,
momenta $p_i\in\C$ and spin vectors and covectors $v_i\in \C^k$,
$u_i\in (\C^k)^*$ ($1\leq i \leq n$), all up to the simultaneous
action of the symmetric group $S_n$. Let $f_{ij}=u_i(v_j)\in\C$
($i\neq j$) be the contraction of the $i$th covector with the $j$th
vector. The hamiltonian for the spin Calogero-Moser system is given
by
$$H= H_2=\frac{1}{2}\sum_{i=1}^n p_i^2 + \sum_{i<j} f_{ij}
f_{ji}U(q_i-q_j).$$
Here the potential function $U$ has a single second
order pole at the origin of $\G$: in terms of coordinates on
the complex line (the universal cover of $\G$), $U$ has one of the forms
\begin{center}
\begin{tabular}{|l|l|}\hline
{\bf Rational:} &  $\displaystyle U(q)=\frac{1}{q^2}$,\\ \hline
{\bf Trigonometric:} & $\displaystyle U(q)=\frac{1}{\sin^2(q)}$,\\ \hline
{\bf Elliptic:}  & $\displaystyle U(q)=\wp(q)$\\ \hline
\end{tabular}
\end{center}
where $\wp(q)$ is the Weierstrass $\wp$-function attached to the
elliptic curve $E$. The {\em spinless} case $k=1$, $f_{ij}=1$ is the
classical Calogero-Moser particle system.

The spin Calogero-Moser hamiltonian depends only on the
contractions $f_{ij}=u_i(v_j)$ (as do the higher integrals of motion
discussed below). Thus the dynamics of the system descend to the
phase space of the {\em Euler-Calogero-Moser system} (we follow
the terminology of \cite{BBKT}), in which we only keep track of the
$p_i,q_i$ and the matrix $F=(f_{ij})\in\gl_n$ (considered up to the
addition of diagonal matrices). This latter system is often referred to
as the spin Calogero-Moser system, since its dynamics come from
those of the full spin CM system. However, the phase spaces for the
two systems are quite different (especially so when $k>n$) and it is
the full spin CM phase space that plays a role in the correspondence
with the multicomponent KP hierarchy \cite{solitons}.

\subsection{CM Matrices}\label{CM matrices}
We briefly recall the description of rational and trigonometric spin
Calogero-Moser systems in terms of matrices (or quivers) following
\cite{KKS,Wilson CM}.

Consider the cotangent bundles $T^*\gl_n^{rss}$ (rational case) and
$T^*GL_n^{rss}$ (trigonometric case) of the regular semisimple loci
in the Lie algebra and group. These cotangent bundles are identified
with the sets of pairs of matrices $(X,Y)$ with $X$ having $n$
distinct eigenvalues ($q_i\in\C$ in the rational/Lie algebra case, $q_i\in\Cx$
in the trigonometric/Lie group case). We now pass to the quotient by the
simultaneous conjugation action of $GL_n$, which is identified with
the phase space of the rational (respectively trigonometric)
Euler-Calogero-Moser system, with $q_i$ the positions of the
particles. The corresponding momenta $p_i$ are recovered as the
diagonal entries of $Y$ in the gauge where $X$ is diagonal. Finally
we have moment maps for the action of $GL_n$: \begin{eqnarray*}
T^*\gl_n^{rss}\to \gl_n,&\hskip.3in&
X,Y\mapsto F=[X,Y],\\
T^*GL_n^{rss}\to \gl_n,&\hskip.3in& X,Y\mapsto F=X\inv Y X-Y,
\end{eqnarray*} giving the spin coordinates $f_{ij}$ as the
off-diagonal entries of the matrix $F\in\gl_n$ (or as $X$-rescaled
versions of the off-diagonal entries of $Y$). Note that this phase
space is Poisson, with symplectic leaves labeled by coadjoint
orbits. The Calogero-Moser hamiltonian is given by the
$GL_n$-invariant function $H^{CM}_2=\frac{1}{2}\on{tr}Y^2$. The
hamiltonians $H^{CM}_i=\frac{1}{i}\on{tr}Y^i$ ($i=1,2,\dots$) are in
involution, and define a degenerately integrable hamiltonian system
\cite{Res}.

 We may now
drop the assumption that the matrix $X$ is regular semisimple,
obtaining a partially completed phase space for the rational and
trigonometric Euler-Calogero-Moser systems in which the positions
$q_i$ are allowed to coincide.

\begin{defn}
The rational and trigonometric Euler-Calogero-Moser spaces\footnote{Here and elsewhere in the paper, we will use the word ``spaces''
(or ``moduli spaces'') in a slightly abusive way.  However, all statements
in the paper apply equally well to the moduli stack
 and to any reasonable moduli spaces/varieties that
result, so the reader may substitute his or her
 preferred type of moduli object.}
 are the
quotients $\ECM_n(\C)=T^*\gl_n/GL_n$, $\ECM_n(\Cx)=T^*GL_n/GL_n$.
The Calogero-Moser hamiltonians on these spaces are the reductions
of the invariant polynomials $H^{CM}_i=\frac{1}{i}\on{tr}Y^i$.
\end{defn}

In order to describe the rational and trigonometric {\em spin}
Calogero-Moser systems, we consider in addition to the matrices
$X,Y$ also maps $u:\C^k\to\C^n$, $v:\C^n\to \C^k$. When $X$ is
regular semisimple, we may write $u,v$ in the basis of $X$-eigenvectors,
giving the data of $n$ vectors $v_i\in \C^k$ and $n$
covectors $u_i\in (\C^k)^*$ that are the spin parameters for the $n$
particles with positions $\{q_i\}$. Noting that the variety of quadruples
$(X,Y,u,v)$ is the cotangent bundle of the variety of pairs $(X,u)$,
we obtain the following definition:

\begin{defn}
The rational and trigonometric spin Calogero-Moser spaces are the
hamiltonian reductions
$$\CM_n^k(\C)=T^*(\gl_n\times\on{Hom}(\C^k,\C^n))/\!\!/_{\on{Id}}GL_n,$$
$$\CM_n^k(\Cx)=T^*(GL_n\times\on{Hom}(\C^k,\C^n))/\!\!/_{\on{Id}}GL_n$$ at
the coadjoint orbit $\{\on{Id}\}\in\gl_n$. In other words, these are
the varieties of quadruples $\{X,Y,u,v\}$ with $[X,Y]+u(v)=\Id$ and
$X\inv Y X-Y+u(v)=\Id$, respectively, modulo the simultaneous action
of $GL_n$. The Calogero-Moser hamiltonians on these spaces are the
reductions of the invariant polynomials
$H^{CM}_i=\frac{1}{i}\on{tr}Y^i$.
\end{defn}
\begin{remark}
It is clear from the above description that
 $\CM_n^k(\C)$ is a (framed) quiver variety \cite{Nakajima} associated to the
quiver with one vertex and one loop (the matrix $X$): the vector $v$
is the framing datum, and $Y,u$ come from doubling the resulting
quiver.
\end{remark}

\subsection{Formulas for Rational CM Matrices} When $k=1$, the spin Calogero-Moser system reduces to
the usual spinless Calogero-Moser system, which is the symplectic
leaf of the Euler-Calogero-Moser system corresponding to the
minimal coadjoint orbit $\Of=\{\Id-u(v)\}\subset\gl_n$ that consists
of traceless matrices of the form ``$\Id$ minus a rank one matrix.''
In other words,
 the (spinless) rational CM space is $$\CM_n= \lbrace (X,Y)\in
T^*\gl_n \;\big|\; [X,Y] \in \Of\rbrace/GL_n.$$
It is proven in
\cite{Wilson CM} that this space is a smooth, irreducible affine variety
of dimension $2n$. It is convenient to realize $\Of$ as the orbit of
the matrix
\begin{equation}\label{CM orbit}\left(\begin{array}{ccccc}
0&1&1&\cdots&1\\ 1&0&1&\cdots&1\\ 1&1&0&\cdots&1\\
\vdots&\vdots&\hdots&\ddots&\vdots\\ 1&1&1&\cdots&0
\end{array}\right) .
\end{equation}
On the open subset where $X$ has distinct eigenvalues, we may then
write coordinates $(q_i,p_i)$ on (a finite cover of) $\CM_n$:
$$X=\left(\begin{array}{ccccc} q_1&0&0&\cdots&0\\ 0&q_2&0&\cdots&0\\
0&0&q_3&\cdots&0\\
\vdots&\vdots&\hdots&\ddots&\vdots\\
0&0&0&\cdots&q_n
\end{array}\right), \hskip.2in
Y=\left(\begin{array}{ccccc}
p_1&\frac{1}{q_1-q_2}&\frac{1}{q_1-q_3}&\cdots&\frac{1}{q_1-q_n}\\
\frac{1}{q_2-q_1}&p_2&\frac{1}{q_2-q_3}&\cdots&\frac{1}{q_2-q_n}\\
\frac{1}{q_3-q_1}&\frac{1}{q_3-q_2}&p_3&\cdots&\frac{1}{q_3-q_n}\\
\vdots&\vdots&\hdots&\ddots&\vdots\\
\frac{1}{q_n-q_1}&\frac{1}{q_n-q_2}&\frac{1}{q_n-q_3}&\cdots&p_n
\end{array}\right).
$$
It is easy to see that the hamiltonian $H=H_2$ in these coordinates
recovers the rational Calogero-Moser hamiltonian above. Thus, $\CM_n$
provides a completion of the phase space of the rational
Calogero-Moser system in which we allow the points $q_i$ (the
eigenvalues of $X$) to collide. In the general spin case, we obtain
coordinates $p_i,q_i$ as diagonal entries just as before and
coordinates $f_{ij}$ ($i\neq j$) with
$Y_{ij}=\frac{f_{ij}}{q_i-q_j}$.

\section{Calogero-Moser and Hitchin Systems}\label{CM and Hitchin chapter}

\subsection{Cubic Curves}\label{cubics}
The classification of $1$-dimensional complex algebraic groups is
paralleled by the classification of Weierstrass cubic curves, that is,
irreducible, reduced complex projective curves $E$ of arithmetic genus $1$:
\begin{enumerate}
\item[$\bullet$] {\bf Elliptic:} $E$ is a smooth elliptic curve (in
particular a group), and may be described by an equation of the form
$y^2=x^3+a x+b$ where $4a^3 + 27b^2\neq 0$.
\item[$\bullet$]{\bf Trigonometric:} $E$ is a nodal cubic, and
is isomorphic to the curve $y^2=x^2(x-1)$.  Its normalization
$\pline\to E$ identifies two points $0$ and $\infty$ to a node on
$E$, and defines a group structure $\Cx= \G\subset E$ on the smooth
locus.
\item[$\bullet$]{\bf Rational:} $E$ is a cuspidal cubic, and is
isomorphic to the curve $y^2=x^3$. Its normalization $\pline\to E$
collapses $2\cdot\infty$ to a cusp on $E$, and defines a group
structure $\C= \G\subset E$ on the smooth locus.
\end{enumerate}

We will denote the identity element of each group $\G$ by $b$,
and the singular point (cusp or node) by $\infty$. In all three
cases, the smooth locus $\G$ is identified (as a group) with the
Jacobian $\on{Pic}^0(E)$ via the map $q\mapsto \cL_q=\Oo(q-b)$. $E$
itself is identified with its compactified Jacobian, the moduli space of
torsion free sheaves of rank $1$ and degree $0$ on $E$.
 The singular point corresponds to the unique rank
1, degree $0$ torsion-free sheaf that is not locally free, namely
the modification ${\mathfrak m}_\infty(b)$ of the ideal sheaf of
$\infty$.

The group variety $\G$ acts on $E$, defining the unique nonzero
invariant vector field $\del$ on $E$ up to a scalar. Writing the
singular cubics in terms of their normalization $\pline$, a choice
of $\del$ is represented by $\frac{\partial}{\partial z}$ in the
cuspidal case (vanishing to order $2$ at $\infty\in\pline$) and by
$z\frac{\partial}{\partial z}$ in the nodal case (vanishing to order
one at $0,\infty$). We will abuse notation to denote the sheaf
$\Oo_E\cdot\del$ by $\cT_E$ (note in the nodal case this is the log
tangent bundle of $E$). Similarly the dual sheaf will be denoted by
$\Omega_E$. Both sheaves are trivial line bundles on $E$.
 The total
space of $\Omega_E$ (which is isomorphic to $E\times \C$) will be
denoted $T^*E$.

\subsection{From Particles to Vector Bundles}\label{particles as
bundles} We would like to encode the positions $q_i\in\G$ of the
Calogero-Moser particles in a ``Fourier dual" fashion.
Recall that $\G$ is identified with the Jacobian $\on{Pic}^0(E)$ of
the corresponding cubic curve $E$, via the map $q\mapsto
\cL_q=\Oo(q-b)$. Thus, $n$ distinct points in $\G$ define a rank $n$
vector bundle $W=\bigoplus \Oo(q_i-b)$ on $E$ which is
semistable of degree zero. Conversely, a generic degree $0$
semistable bundle on $E$ is of the form $W=\bigoplus \Oo(q_i-b)$ for
$n$ distinct points $q_i\in\G\subset E$ (determined up to
permutation).

On an elliptic curve $E$, we may extend this correspondence as
follows. The Fourier-Mukai transform \cite{Mukai} identifies
semistable degree zero vector bundles $W$ with degree (length) $n$ torsion
coherent sheaves $W^\vee$ on $E$, in such a way that $W=\bigoplus
\Oo(q_i-b)$ is identified with the sum of skyscrapers $\bigoplus
\Oo_{q_i}$. We may therefore consider the moduli space of such $W$
as a completion of the configuration space of $n$ points $q_i\in E$.

A similar argument for singular cubic curves identifies certain
semistable degree zero vector bundles with configurations of points
in the {\em smooth} locus $\G\subset E$:

\begin{lemma}[\cite{FriedMorgminuscule}, Lemma 1.2.5]\label{trivial pullbacks}
Suppose $W$ is a semistable vector bundle of degree $0$ on a
singular Weierstrass cubic $E$, and that $\pline\simeq
\widetilde{E}\xrightarrow{n} E$ is the normalization. Then the
following are equivalent:
\begin{enumerate}

\item $W$ is identified by the Fourier-Mukai transform with a torsion coherent
sheaf $W^\vee$ supported on $\G\subset E$.

\item $n^*W$ is a trivial vector bundle.

\item $W$ has a filtration whose subquotients are line bundles of degree
$0$ on $E$.
\end{enumerate}
\end{lemma}

For any cubic curve $E$, we let $Bun_n^{ss}(E)$ denote the moduli
stack of rank $n$ semistable vector bundles of degree zero on a
cubic curve $E$; for singular $E$, we impose in addition any of the
equivalent (and open) conditions of Lemma \ref{trivial pullbacks}.
Note that we do not at any point pass to the moduli {\em space} of
($S$-equivalence classes of) semistable vector bundles.

\subsection{Calogero-Moser and Hitchin Systems}\label{CM and
Hitchin}

In this section we review the identification of the completed
Euler-Calogero-Moser systems with meromorphic Hitchin systems on
cubic curves following \cite{Nekrasov} (see also
\cite{GorNe,Enriquez,BBT}).

We denote by $Bun_n^{ss}(E,b)\to Bun_n^{ss}(E)$ the principal
$GL_n$-bundle parametrizing bundles in $Bun_n^{ss}(E)$ equipped with
a trivialization of the fiber at the identity.  The cotangent fiber
$T^*Bun_n^{ss}(E,b)|_W$ at a bundle $W$ consists of pairs
$(W,\eta)$, where $\eta$ is a meromorphic Higgs field
$\eta\in\Gamma(\on{End}W(b))$ on $W$ with only a simple pole at $b$.
The group $GL_n$ admits a hamiltonian action on $T^*Bun_n^{ss}(E,b)$
induced from its action on $Bun_n^{ss}(E,b)$ (by changing the
trivialization at $b$). Let $\C_{k,n}=\on{Hom}(\C^k,\C^n)$, and
identify $T^*\C_{k,n}=\C_{k,n}\times \C_{n,k}$.

\begin{defn}\label{CM spaces}
The Hitchin-Calogero-Moser space $\HCM_n(E)$ associated to the
cubic curve $E$ is the quotient $T^*Bun_n^{ss}(E,b)/GL_n$, i.e. the space
of Higgs bundles on $E$ with simple pole at $b$.
\end{defn}

\begin{defn}\label{framed Higgs} A {\em framed Higgs bundle} is a
quadruple $(W,\eta,u,v)$ where
\begin{enumerate}
\item $(W,\eta)\in \HCM_n(E)$ is a Higgs
bundle on $E$ with pole at $b$,
\item  $u:\C^k\to W|_b$ and $v:W|_b\to \C^k$ are linear maps,
\end{enumerate}
and we require that $\on{Res}_b\eta+u(v)=\Id$.

The moduli space of framed Higgs bundles on $E$ is denoted
$\HCM_n^k(\G)$ and is identified with the hamiltonian reduction
$$\HCM_n^k(\G)=T^*(Bun_n^{ss}(E,b)\times\C_{k,n})/\!\!/_{\Id} GL_n.$$
\end{defn}

As we explain below, the Hitchin systems on the spaces $\HCM_n(E)$
give a completion of the Euler-Calogero-Moser systems. The spaces
$\HCM_n^k(E)$ of framed Higgs bundles, on the other hand, will model
the spin CM systems, and the map $\HCM_n^k(E)\to \HCM_n(E)$
forgetting $u,v$ corresponds to forgetting the spins and remembering
only their contractions $f_{ij}$. In particular the Hitchin
hamiltonians on $\HCM_n(E)$ pull back to define the framed Hitchin
system on $\HCM_n^k(E)$. In Section \ref{spectral sheaves and higgs}
we will relate framed Higgs bundles (and their generalizations) to
spectral sheaves equipped with a normalization
(framing) on a ruled surface.

\begin{remark}
The contraction map $\HCM_n^k(E)\to \HCM_n(E)$ forgetting $(u,v)$
identifies the (spinless) Calogero-Moser space
$\CM_n(E)=\HCM_n^1(E)$ with the subspace (in fact symplectic leaf)
of $\HCM_n(E)$ consisting of Higgs fields with residue in the
coadjoint orbit $\Of$.
\end{remark}

\subsection{Matrices and Hitchin Systems}\label{matrices} The spaces of Calogero-Moser
matrices $\ECM_n(\C)$ and $\ECM_n(\Cx)$ are readily identified with
the rational and trigonometric Hitchin-Calogero-Moser spaces
$\HCM_n(E)$. Namely, by hypothesis, bundles in $Bun_n^{ss}(E,b)$ have
trivial pullback to the normalization $\pline$, so are completely
described by the descent data from $\pline$ to $E$. This descent
datum in the nodal case is the identification of the two fibers over
the inverse image of the node, hence $Bun_n^{ss}(E,b)=\GL_n$. In the
cuspidal case, these two points are infinitesimally nearby, and the
descent datum becomes a ``connection matrix'' identifying these two
nearby fibers---thus, we have $Bun_n^{ss}(E,b)=\gl_n$. (See
\cite{FriedMorgIII} for more details.)

In general, fix a cubic curve $E$ and an invariant differential on
it. We restrict to the open locus in $Bun_n^{ss}(E)$ consisting of
vector bundles $W\simeq \bigoplus\cL_{q_i}$, sums of the line
bundles $\cL_{q_i}=\Oo(q_i-b)$ of degree zero associated to $n$
distinct points $q_i\in \G$. Let $s_{q_i-q_j}$ denote the unique
section meromorphic section of $\cL_{q_i-q_j}$ with only a simple
pole at $b$ (normalized using the differential and trivialization of
the fiber) and zero at $q_i-q_j$. Then it is easy to see that the
Higgs field $\eta\in \End(W)(b)=\bigoplus\cL_{q_i-q_j}(b)$ must have
the form
\begin{equation}\label{form of Higgs}
\left(
\begin{array}{ccccc}
p_1&f_{12}s_{q_1-q_2} & f_{13} s_{q_1-q_3}&\cdots&f_{1n} s_{q_1-q_n}\\
f_{21}s_{q_2-q_1}&p_2&f_{23}s_{q_2-q_3}&\cdots&f_{2n}s_{q_2-q_n}\\
f_{31}s_{q_3-q_1}& f_{32}s_{q_3-q_2}&p_3&\cdots&f_{3n}s_{q_3-q_n}\\
\vdots&\vdots&\hdots&\ddots&\vdots\\
f_{n1}s_{q_n-q_1}&f_{n2}s_{q_n-q_2}&f_{n3}s_{q_n-q_3}&\cdots&p_n
\end{array}
\right)\end{equation} for some $f_{ij}\in \C$ and where the $p_i$
are sections of $\Oo(b)$ on $E$, hence constants. In the rational
case we have $s_{q_i-q_j}=\frac{1}{z}-\frac{1}{q_i-q_j}$, so that
writing $X=diag(q_i)$ and $[X,Y]=(f_{ij})$ we have
\begin{equation}\label{rational Higgs}
\eta=(\frac{[X,Y]}{z}+ Y)dz,\end{equation} with $Y$ as in equation
\ref{form of Higgs} with $s_{q_i-q_j}$ replaced by
$\frac{1}{q_i-q_j}$. Similarly in the trigonometric case we replace
$\frac{1}{q_i-q_j}$ by $\sin(q_i-q_j)$ and
\begin{equation}\label{trig Higgs}
\eta=(\frac{X\inv YX-Y}{z}+ Y)\frac{dz}{z}.
\end{equation}

\subsection{CM Hamiltonians} Finally, we give the rational,
trigonometric and elliptic spin Calogero-Moser hamiltonians in the
Hitchin system description. Recall that the Hitchin hamiltonians on
the moduli space of Higgs bundles are all components of traces of
powers of the Higgs field. Among these we wish to single out the
spin CM hamiltonians.

\begin{defn}\label{spin hamiltonians} The spin CM hamiltonians on $\HCM_n(E)$ are the
functions $$H_i:(W,\eta)\mapsto\frac{1}{(i+1)}\on{Res}_b
\on{Tr}(\eta+\zeta\on{I})^{i+1}.$$
\end{defn}

We see from the explicit form of the Higgs fields in the rational
and trigonometric cases (Equations \ref{rational Higgs} and
\ref{trig Higgs}) that
$\frac{1}{i}H_i(X,Y)=H_i^{CM}(Y)=\frac{1}{i}\on{tr}Y^i$ are indeed
the CM hamiltonians.\footnote{Our normalization of the hamiltonians
is chosen to be compatible with the Hitchin hamiltonians in the next
section.} Similarly, in the elliptic case one checks that
$\frac{1}{2}H_2$ is the (quadratic) elliptic CM hamiltonian using
the identity $s_{q_i-q_j}s_{q_j-q_i}=\wp$.

Summarizing, we have the following statement.

\begin{prop}[\cite{Nekrasov}]\label{nekrasov} The
Hitchin-Calogero-Moser systems on cuspidal, nodal, and smooth
cubic curves are completions of the rational, trigonometric and
elliptic  Euler-Calogero-Moser systems, respectively. Moreover, we
have isomorphisms $\ECM_n(\C)\simeq \HCM_n(\C)$ and
$\ECM_n(\Cx)\simeq \HCM_n(\Cx)$ of integrable systems.
\end{prop}

\section{Framed Calogero-Moser Systems and Spectral Sheaves}\label{CM spectral
sheaves}

\subsection{The Surface $\Enat$ and the Weierstrass
$\zeta$-Function.}\label{about Enat}
In
this section, we discuss a surface $\Enat$ that is the total space of
the unique (up to
isomorphism)
nontrivial rank one affine bundle over a cubic curve $E$.
\subsubsection{$\Enat$ for Smooth Cubics}
 Fix an elliptic curve $E$.  Let
$\At$ denote the {\em Atiyah bundle} on $E$, the unique
(again, up to isomorphism) nontrivial
extension of $\Oo_E$ by itself. The
algebraic surface $\Enat$ is the complement of the
section $E_\infty=\Pp(\Oo)\cong E$ of the projectivization $\Enatbar = {\mathbb P}(A)$ of the
Atiyah bundle:
\bd
\Enat= \Enatbar \setminus E_\infty.
\ed
 The
resulting surface $\Enat$ is Stein (but not affine algebraic),
and isomorphic with the moduli
space of line bundles with a holomorphic connection on $E$ (see
\cite{TV,TV2} for more on the geometry of $\Enat$).

In classical analytic terms, the surface $\Enat$ may be viewed as
the receptacle for the Weierstrass $\zeta$-function of $E$, the
unique odd function on the universal cover $\C$ of $E$ whose
derivative $\zeta'(z)=-\wp(z)$ is minus the Weierstrass
$\wp$-function of $E$. That is, while $\zeta$ is not
doubly-periodic, it differs from its translates by additive
constants, so that it determines a well-defined section of an affine
$\C$-bundle over $E$. This surface is readily identified with
$\Enat$. Indeed,  recall that the Weierstrass $\sigma$-function of
$E$ is the unique section of the line bundle $\Oo(b)$ with a simple
zero at $b$, and that $\zeta=d\on{log}\, \sigma(z)$ is its
logarithmic derivative. This provides the algebraic definition of
$\zeta$: it is the section of the affine bundle ${\mc Conn} \Oo(b)$
of connections on $\Oo(b)$ with log pole at $b$ and residue $1$ that
corresponds to the unique meromorphic connection annihilating
$\sigma$.

We may now fix $\Enat$ up to unique isomorphism by setting
$\Enat={\mc Conn}\Oo(b)$, a {\em twisted cotangent bundle} of $E$
\cite{BB} (i.e. affine bundle for $\Omega_E\simeq \Oo_E$ with
compatible symplectic structure). Let $\cA$ denote the pushforward
of $\Oo_{\Enat}$ to $E$, i.e. the algebra of functions on the fibers
of $E$. Thus $\At=(\cA)_{\leq 1}$, the subsheaf of affine functions
on $\Enat$, is isomorphic to the Atiyah bundle, and is canonically
an extension of $\cT_E$ (which is isomorphic to $\Oo_E$) by $\Oo_E$.
The sheaf $\At$ is also isomorphic as $\Oo_E$-module to $\D^1(\Oo_E(b))$,
the sheaf of differential operators of order at most
one acting on the line bundle $\Oo(b)$. Concretely, the sheaf $\At$
lies in between
\begin{equation}\label{bracketing At}
\Oo_E\oplus\cT_E(-b)\subset \At \subset \Oo_E\oplus \cT_E(b),
\end{equation} and $\At$ is generated (in the canonical local coordinate
near $b$) by $\Oo_E\oplus\cT_E(-b)$ and the section $\del-\zeta$.

The meromorphic section $\zeta(z)$ of $\Enat\to E$ defines a
trivialization of the affine bundle $\Enat\to E$ away from $b$, and
hence a canonical birational identification of the cotangent bundle
$T^*E\simeq E\times\C$ with $\Enat$:
 $$T^*E\ni
(z,\omega(z))\mapsto(z,\zeta(z)+\omega(z))\in \Enat$$ (see \cite{DM}
for a geometric description in terms of elementary modifications).
Let $\ul{k}$ denote the meromorphic function on $\Enat$, with polar
divisor the fiber over $b$, obtained by composing this
identification with the projection onto $\C$,
$\ul{k}(z,\zeta(z)+\omega(z))=\omega(z)$ (or more canonically
$\omega(z)/dz$). Let $\underline{\zeta}$ denote the Laurent
expansion of $\zeta$ at the origin $0\in\C$, considered as a Laurent
series on $E$ at $b$. Then it follows that the function
$\underline{t}=\ul{k}+\pi^*\underline{\zeta}$
on $\Enat$ near $F_b$ is regular along the fiber $F_b$, and gives a
natural affine coordinate on $\Enat$ near $F_b$ (i.e an affine
identification of the formal neighborhood of $F_b\subset \Enat$ with
$E\times \C$). More generally, pullback by $\ul{k}$ identifies
regular functions on $\Enat$ with the ring of functions generated by
linear functions $f$ on $T^*X$ such that $f+\pi^*\ul{\zeta}$ is
regular near $F_b$ (as is evident from equation \ref{bracketing
At}).

\subsubsection{$\Enat$ for Singular Cubics}
The definition and properties of $\Enat$ extend naturally to general
cubics $E$. For any Weierstrass cubic $E$ we have \bd
\Ext^1_E(\theo,\theo)= H^1(E,\theo) \cong \C. \ed So $E$ has a
unique nontrivial extension $\At$ of $\theo_E\simeq \cT_E$ by
$\theo_E$, up to isomorphism.\footnote{Recall that $\cT_E$ is the subsheaf
of the tangent sheaf generated by the $\G$-action.} We again fix
$\At=\D^1(\Oo_E(b))$ (these are differential operators with
symbol in $\cT_E$). Let $\Enatbar \overset{\on{def}}{=} \bproj(\Sym
\At)$ denote the associated ruled surface, and $p:\Enatbar
\rightarrow E$ denote the projection map.

The quotient map $\At\twoheadrightarrow \theo_E$ defines a section
$s: E\rightarrow \Enatbar$; we write $E_\infty = s(E)$ and refer to
it as the {\em section at infinity}.  The surface $\Enat
\overset{\on{def}}{=} \Enatbar\smallsetminus E_\infty$ is called the
{\em twisted (log) cotangent bundle} of $E$; it is the nontrivial
torsor over $\Omega_E$ given by the nonzero class (up to scale) in
$H^1(\Omega_E)$.

The notion of $\zeta$-function and its relation to $\Enat$ similarly
extend to the singular cases (again considering $\zeta$ as a log
connection on $\Oo(b)$). Concretely, in the rational case, we have
$\zeta=1/z$ and in the trigonometric case
$\zeta=1/\sin(z)$.

\subsection{Twisted Higgs Fields and Framed CM Systems}\label{twisted Higgs}
In this section we introduce the notion of {\em twisted} Higgs
fields, which are a modified version of Higgs fields whose spectral
curves naturally live in $\Enat$ rather than in the cotangent
bundle. The two notions are readily identified, but the translation
from twisted Higgs to Higgs adds $\on{Id}$ to the residue at the
basepoint, providing a geometric origin to the appearance of
$\zeta\on{Id}$ in the CM Hamiltonian or equivalently of $\on{Id}$ in
the CM moment condition. We then define the framed CM particle
system in its Hitchin system formulation.

Recall (Section \ref{about Enat}) the $\Oo_E$-algebra $\cA$ of
functions on $\Enat$ and its subsheaf $\At$ (the Atiyah bundle) of
affine functions, which is an extension of $\Oo_E$ by $\cT_E\simeq
\Oo_E$. The structure of Higgs bundle on a vector bundle $W$ can be
written as an action of $\Oo_E\oplus \cT_E$ on $W$ extending the
$\Oo_E$-module structure, which makes $W$ into a sheaf on $T^*E$.
This has a natural twisted analog, in which we replace the
$\Oo_E\oplus\cT_E$-action by an action of $\At$, or equivalently a
lifting to a sheaf on $\Enat$:

\begin{defn} A (regular) {\em twisted Higgs bundle} on $E$ is a pair
$(W,\wt{\eta})$ where 
\begin{enumerate}
\item $W$ is a vector bundle on $E$ and
\item $\wt{\eta}:
\At\ot_{\Oo_E} W\to W$ is a map whose restriction to $\Oo_E\subset
\At$ is the identity map of $W$.
\end{enumerate}
\end{defn}

The relation between Higgs and twisted Higgs fields is given by the
Weierstrass $\zeta$-function, which may be considered as a splitting
of the extension $\At$ away from $b$, or as a section of $\At$ with
simple pole at $b$. Put another way, $\At$ is identified with the
subsheaf of $\Oo_E\oplus\cT_E(b)$ generated by $\Oo_E$, $T_E(-b)$
and the section $\del-\zeta$. It follows that to give a twisted
Higgs field $\wt{\eta}$ on $W$ is equivalent to giving an action of
$\del$ from $W$ to $W(b)$, which becomes regular after subtracting
$\zeta\Id$. In coordinate-free language, this is a Higgs field
$\eta=\wt{\eta}-\zeta\on{I}$ on $W$ with simple pole whose residue
at $b$ is the identity endomorphism. Recall (Section \ref{about
Enat}) that we have a canonical birational isomorphism between
$T^*E$ and $\Enat$ relative to $E$, given by the function $\ul{k}$.

\begin{lemma}
Let $W$ denote a vector bundle on $E$. There is a
bijection between meromorphic Higgs fields $\eta:W\to W(*b)$ and
meromorphic twisted Higgs fields $\wt{\eta}:\At\ot_{\Oo_E} W\to
W(*b)$ sending $\eta$ to $\wt{\eta}-\zeta\on{I}$. The corresponding
spectral sheaves on the surfaces $T^*E$, $\Enat$ away from the fiber
over $b$ are identified by the birational isomorphism of the
surfaces given by $\zeta$.
\end{lemma}

\subsubsection{Framed Higgs Fields}
Let $\framing$ denote a torsion coherent sheaf on $E$, with support
$S\subset \G$ a subscheme of the smooth locus of $E$ (considered as
a divisor on $\G$). For a sheaf $F$ on $E$ we use the chosen
invariant differential on $E$ to identify $F(S)/F$ with the
restriction $F|_S$.

\begin{defn}\label{framed twisted Higgs def}
 A $\framing$-{\em framed twisted Higgs bundle} is
a quadruple $(W,\wt{\eta},u,v)$ where
\begin{enumerate}
\item $W\in Bun_n^{ss}(E)$,
\item 
$u:\framing\to W|_S$ and  $v:W|_S\to \framing$ are maps of coherent sheaves, 
and
\item $\wt{\eta}:\At\ot_{\Oo_E} W\to W(S)$ is a map whose restriction to
$\Oo_E\subset \At$ is the identity map of $W$;
\end{enumerate}
these data must satisfy the following.   The restriction of
$\wt{\eta}$ to $S$ factors through a map
$$p.p.(\wt{\eta}):W|_S = (A\otimes W/\theo\otimes W)|_S\to W(S)|_S\simeq W|_S$$ which we require to satisfy
\bd
p.p.(\wt{\eta})=u(v).
\ed
 We denote the moduli space of framed twisted
Higgs bundles by $\HCM_n(E,\framing)$.
\end{defn}

\begin{lemma}\label{twisted CM}
There is an isomorphism between the moduli spaces of:
\begin{enumerate}
\item $\framing$-framed twisted Higgs bundles $(W,\wt{\eta},u,v)$, and
\item quadruples $(W, \eta, u, v)$ where
\begin{enumerate}
\item
$\left(W\in
Bun_n^{ss}(E) , \eta:W\to W(S+b)\right)$
is a meromorphic Higgs bundle and
\item  $u:\framing\to
W|_S$ and $v:W|_S\to \framing$ are maps of coherent sheaves
\end{enumerate}
 whose principal parts satisfy
$$p.p.(\eta)=u(v)+p.p.(\zeta\cdot\Id_W):W|_S\to W|_{S+b}.$$
\end{enumerate}
\end{lemma}

In the case $\framing=\Oo_b^k$, we recover the notion of framed
Higgs bundle from Definition \ref{framed Higgs}: $$\HCM_n^k(E)\simeq
\HCM_n(E,\Oo_b^k).$$ The notion of twisted Higgs bundle thus gives a
geometric interpretation (the passage from $T^*E$ to $\Enat$) for
the appearance of $\Id$ in the CM residue condition on Higgs fields.

\subsubsection{Framed CM Systems}\label{framed CM systems}
We now define hamiltionians on the framed Calogero-Moser phase
spaces $\HCM_n(E,\framing)$ as certain Hitchin hamiltonians,
generalizing Definition \ref{spin hamiltonians} in the spin case
$\framing=\Oo_b^k$. Let $S=\coprod_{j=1}^\ell x_j\subset E$ denote the
set-theoretic support of $\framing$, consisting of $\ell$ distinct
points on $X$ (counted without multiplicity). We define hamiltonians
$H_{i,b}$ whenever $b\in S$ as in the spin case, but also a
collection of hamiltonians $H_{i,x_j}$ for each point $x_j\in S$, as
follows:

\begin{defn}\label{framed CM hamiltonians}
The framed CM hamiltonians are the functions
$$H_{i,x_j}:\HCM_n(E,\framing)\to \C, \hskip.3in
H_{i,x_j}(W,\eta,u,v)=\frac{1}{i+1}\on{Res}_{x_j}\on{Tr}
\eta^{i+1}$$ for $x_j\neq b$, together with (when $b\in S$)
$$H_{i,b}(W,\eta,u,v)=\frac{1}{i+1}\on{Res}_{b}\on{Tr}(\eta+\zeta\on{I})^{i+1}.$$
\end{defn}

The framed CM systems may be identified with completed particle
systems on the group $\G$, as in Section \ref{matrices}. Namely, we
restrict to the open locus in which the vector bundle $W=\bigoplus
\Oo(q_i-b)$ is canonically a sum of distinct line bundles (up to
permutation). The positions of the corresponding particles are then
given by the points $q_i\in \G$ (given by the Fourier-Mukai
transform of $W$, as in Section \ref{particles as bundles}). The
decomposition of $W$ moreover allows us to decompose $\eta$ into
diagonal and off-diagonal components. We identify the momentum $p_i$
of the particle $q_i$ as the constant term of the $i$th diagonal
component of $\eta$. It is then immediate that the hamiltonian
$\frac{1}{2}H_{2,b}$ consists of the kinetic term $\frac{1}{2}\sum
p_i^2$ together with other (potential) terms. The other hamiltonians
above define integrals of motion for this system. In the case of
{\em simple framing} (that is, when $\framing=\bigoplus x_i$ for distinct
$x_i$) they form a maximal family of integrals in involution,
defining an algebraically completely integrable system. (This
follows immediately from the spectral description of framed CM
systems, Theorem \ref{framed CM flows thm}.)

\subsection{CM Spectral Sheaves}\label{CM spectral sheaves section}
Moduli spaces of spectral sheaves (specifically, of line bundles on
curves in a Poisson surface) give a wide class of examples of
integrable systems (see e.g. \cite{DM,Hu}). The prototypical example
of such a setting is the ($GL_n$) Hitchin system on the moduli space
$T^*Bun_n(X)$ of Higgs bundles on a curve $X$, which can be
described as a moduli space of torsion-free sheaves on curves in $T^*X$
finite of degree $n$ over $X$. We will similarly realize the spin CM
systems in terms of spectral curves on the twisted cotangent bundle
$\Enat$ of $E$. More precisely, we look at torsion-free sheaves on
curves in the projective surfaces
\bd
\ol{T^*E}=T^*E \cup E_\infty \;\;\;\text{and}\;\;\; \Enatbar=\Enat\cup
E_\infty,
\ed
 for which we fix the behavior along the curves
$E_\infty\equiv E$ at infinity.

\begin{defn}\label{CM spectral sheaves def} Fix a
coherent torsion sheaf $\framing$ on $\G$. A $\framing$-{\em framed CM
spectral sheaf} (respectively, {\em Hitchin spectral sheaf}) is a pair
$(\cF, \phi)$ consisting of a coherent sheaf $\cF$ on $\Enatbar$
(respectively, $\ol{T^*E}$) of pure dimension one, together with an
identification $\phi:\cF|_{E_{\infty}}\to \framing$, satisfying the
following two normalization conditions:
\begin{enumerate}
\item[(i)] $W=p_*\cF(-E_\infty)$ is a semistable vector bundle of
degree $0$; if $E$ is singular, we also require that the pullback of
$W$ to the normalization of $E$ is a trivial vector bundle.
\item[(ii)] $\on{deg}(p_*\cF(kE_\infty))=(k+1)\on{deg}(\framing)$ for
$k\gg 0$.
\end{enumerate}
We denote the sheaf $p_*\cF(kE_\infty)$ by $F_k$. The {\em
$\framing$-Calogero-Moser space $\CM_n(E,\framing)$} is the moduli
scheme of $\framing$-framed CM spectral sheaves $(\cF,\phi)$ for
which the rank of the vector bundle $W$ is $n$.
\end{defn}

As we'll see (Definition \ref{framed twisted Higgs def}, Lemma \ref{twisted CM}, and
Corollary \ref{CM is CM}), the {\em spins}
of generalized Calogero-Moser particles take value in the sheaf
$\framing$. We will identify the moduli space
$\CM_n^k(E):=\CM_n(E,\Oo_b^k)$ with the completed phase space of the
usual $k$-spin $n$-particle Calogero-Moser system. For general
framing, $\CM_n(E,\framing)$ is identified with the framed Hitchin
space $\HCM_n(E,\framing)$. At the other extreme from the spin CM
case $\framing=\Oo_b^k$ we have the case of {\em simple framing},
$\framing=\bigoplus_{i=1}^k \Oo_{x_i}$ with the $x_i$ all distinct,
which also generalizes the spinless case $\framing=\Oo_b$.

\begin{remark}[Normalization Conditions]\label{normalization
conditions} The normalization conditions (i) and (ii) are open
conditions on coherent sheaves of pure dimension one (and in fact
any $\framing$-framed $\F$ satisfying (i) must have pure dimension
one). Condition (i) on the vector bundle $W$ is
discussed in Section \ref{particles as bundles}: such $W$ encode
(via the Fourier transform) the positions of the Calogero-Moser
particles. Condition (ii) is a normalization on the Hilbert
polynomial of $\cF$ and should be considered part of the framing
data. Note that (i) contains the case $k=-1$ of (ii). In fact, it is
easy to see that (i) and (ii) together imply Condition (ii) for all
$k$: we use the exact sequence
$$0 \to F_{k-1} \to F_k \to \framing$$ to conclude that
$\on{deg}(F_k)\leq \on{deg}(F_{k-1}) + \on{deg}(\framing)$.
Then induction on k gives that $\on{deg}(F_k) \leq
(k+1)\on{deg}(\framing)$ provided $\on{deg}(F_{-1})=0$, with
equality if and only if the above sequences are right exact for all
$k\geq 0$.

It is also worth noting that if $\cF$ satisfies hypothesis (i) of
the definition, then it satisfies hypothesis (ii) if and only if the
natural map
\bd
p_*\cF(kE_\infty)\rightarrow
p_*(\cF(kE_\infty)|_{E_\infty})\cong \framing
\ed
 is surjective for all
$k\geq 0$.
\end{remark}

\begin{remark}[Hamiltonians]
The space of framed
CM hamiltonians admit a natural geometric description on the
moduli space of $T$-framed CM spectral sheaves.  Indeed, following
\cite{DM, Hu, T matrix, TV, TV2}, one may take a spectral sheaf to its scheme-theoretic support,
viewed as a divisor on $\Enatbar$ in the linear series
\bd
{\mathbb P}= |\on{rk}(W)\cdot E_\infty + \on{deg}(T)\cdot F|,
\ed
 where $F$ is a fiber of the projection
$\Enatbar\rightarrow E$.  The collection of divisors that contain the curve
$E_\infty$ with nonzero multiplicity form a hyperplane in ${\mathbb P}$ with
complement an affine space ${\mathbb A}$, and we obtain a natural map
$\CM_n(E,\framing)\rightarrow {\mathbb A}$.   The framed CM Hamiltonians
then come by pulling back a particular list of polynomials from ${\mathbb A}$.
\end{remark}

\subsection{From Spectral Sheaves to Higgs
Bundles}\label{spectral sheaves and higgs} In this section we
identify the moduli spaces of framed spectral sheaves on $\Enatbar$
and framed twisted Higgs bundles. This result is based on a ``Koszul
dual" description of framed sheaves due to L. Katzarkov, D. Orlov
and T. Pantev \cite{KOP}.

\begin{thm}[\cite{KOP}]\label{framing theorem}
There is a canonical equivalence between the category of
$\framing$-framed CM spectral sheaves and that of {\em Koszul
data}: quintuples $(W, W', \iota, s, a)$ consisting of
\begin{enumerate}
\item $W\in Bun_n^{ss}(E)$,
\item an extension
$$ 0\to W \xrightarrow{\iota}W' \xrightarrow{s} \framing\to 0$$
of $W$ by $\framing$, and
\item
a map $a:\At\ot W\to W'$ extending $\iota$ on $W=\Oo_E\ot W\subset
\At\ot W$.
\end{enumerate}
\end{thm}
\begin{construction}
The assignment of Koszul data to a spectral sheaf proceeds as follows.
Let $(\cF,\phi)$ be a $\framing$-framed CM spectral sheaf.  We set 
$W= p_*\cF(-E_\infty)$ and $W' = p_*\cF$.  
Then $W$ satisfies Condition (i) of
 Definition
\ref{CM spectral sheaves def}.  Moreover, by Remark \ref{normalization conditions}, the natural sequence
\bd
0\rightarrow W\xrightarrow{\iota} W'\xrightarrow{s} T \rightarrow 0
\ed
is exact.
Making the identification $A = p_*\theo(E_\infty)$, we let $a:A\ot
W\to W'$ be the restriction of the action of $\cA$ on sections of
$\cF$.
\end{construction}

\begin{thm}\label{CM identification}
 There is an isomorphism $\CM_n(E,\framing) \to
\HCM_n(E,\framing)$ between the moduli of $\framing$-framed
spectral sheaves and (untwisted) Higgs fields.
\end{thm}

\begin{proof}
We need to establish an equivalence between
$\framing$-framed twisted Higgs bundles $(W,\wt{\eta},u,v)$ and
Koszul data $(W, W', \iota, s, a)$ as in Theorem \ref{framing
theorem}.

We first establish a bijection between the two types of
 data $(W, W', \iota, s)$ and 
$(W, u)$ coming from Koszul data and Higgs data, respectively. 
Given $(W,W',\iota, s)$, 
there is a natural map $\wt{u}$,
\bd
\xymatrix{
W'\ar@/^1.5em/[rr]^{\wt{u}}\ar@{>>}[r] & W'/\on{tors}(W')\ar[r] & W(S),}
\ed
 induced by the 
isomorphism $W|_{E\setminus S} \cong W'|_{E\setminus S}$.  
The map $\wt{u}$
 restricts to
 the
identity on $W$, and we denote the associated quotient map
$T \xrightarrow{s^{-1}} W'/W \xrightarrow{\wt{u}} W(S)/W$
by $u$; we may also use the canonical identification $W(S)/W = W/W(-S) = W|_S$ coming from the invariant differential of $E$ to 
identify $u$ with a map $T\rightarrow W/W(-S) = W|_S$.

Conversely, a diagram chase shows that $W'$ is obtained (up to
unique isomorphism) as the pullback of the exact sequence
\bd
0\rightarrow W\rightarrow W(S) \rightarrow W(S)/W\rightarrow 0
\ed
along the map $u: T \rightarrow W(S)/W$.  It is immediate that these two 
constructions give the bijections
\begin{equation}\label{W bijections}
\lbrace (W,W',\iota,s)\rbrace \;\;\leftrightharpoons\;\; \lbrace (W,u)\rbrace.
\end{equation}

We thus obtain a diagram
\begin{equation}\label{main Higgs diagram}
\xymatrix{
0\ar[r] & W \ar[r]^{\iota} \ar@{=}[d] & W'\ar@{}[dr]|{\Box} \ar[r]^{s} \ar[d] & T \ar[r]\ar[d]_{u} & 0\\
0 \ar[r] & W \ar[r] & W(S) \ar[r] & W(S)/W \ar[r] & 0,}
\end{equation}
where the square marked $\Box$ is Cartesian, relating the corresponding data
$(W, W', \iota, s)$ and $(W,u)$.  

It is now immediate from the universal propety of pullbacks applied to 
\eqref{main Higgs diagram}
 that there is a bijection between:
\begin{enumerate}
\item the set of 
maps $a: A\otimes W\rightarrow W'$ such that $a|_{\theo\otimes W}$ is the identity on $W$.
\item The set of pairs 
\bd
\left(A\otimes W\xrightarrow{\wt{\eta}} W(S), A\otimes W\xrightarrow{\wt{v}} T\right)
\ed 
such that
\begin{enumerate}
\item
$\wt{\eta}|_{\theo\otimes W}$ is the identity on $W$, 
\item $\wt{v}|_{\theo\otimes W} = 0$, and 
\item the diagram
\bd
\xymatrix{A\otimes W \ar[d]^{\wt{\eta}} \ar[r]^{\wt{v}} & T\ar[d]^{u} \\
W(S)\ar[r] & W(S)/W}
\ed
commutes.
\end{enumerate}
\end{enumerate}
Since the maps $\wt{v}$ in (2) are completely determined by the
induced maps 
\bd
W|_S = (A\otimes W/\theo\otimes W)|_S \xrightarrow{v} T,
\ed
 we find that the bijections of \eqref{W bijections} extend to bijections
\bd
\lbrace (W,W',\iota,s,a)\rbrace \;\;\leftrightharpoons\;\; \lbrace (W,\wt{\eta},u, v)\rbrace,
\ed
as desired.  The proof of the functoriality properties of this bijection necessary to obtain a moduli isomorphism is straightforward, and we omit it.
\end{proof}

\begin{corollary}\label{CM is CM}
The moduli space of $\theo_b^k$-framed spectral sheaves is a completed
phase space for the $k$-spin Calogero-Moser system.
\end{corollary}
\begin{proof}
This is immediate from Theorem \ref{CM identification} and Proposition
\ref{nekrasov}.
\end{proof}

\section{Flows on Spectral Sheaves}\label{flows section}

\subsection{Tweaking Sheaves}\label{tweaking section}
In this section we consider some variants of the simplest method of
deforming sheaves on any variety $Y$, namely tensoring them by line
bundles. If $Y$ is a smooth projective variety, this gives rise to an
action of the Picard group of $Y$ on moduli spaces of sheaves on
$Y$. In particular, the tangent space $\HH^1(Y,\Oo)$ to $\on{Pic}Y$
at the trivial bundle gives rise to infinitesimal deformations of
any sheaf. This infinitesimal action is defined for an arbitrary
variety $Y$ as the canonical map
$$\HH^1(Y,\Oo)=\Ext^1(\Oo,\Oo)\to \Ext^1(\cF,\cF)$$ for any sheaf
$\cF$. Concretely, a self-ext or first-order deformation of $\Oo$
defines, via tensor product, a first-order deformation of any sheaf
$\cF$.

To construct particular deformations of sheaves (on reasonable
varieties $Y$), we can produce elements of $\HH^1(Y,\Oo)$ from local
$\HH^1$ of $\Oo$ along divisors, or from arbitrary meromorphic
functions on $Y$. If $\cF$ is a vector bundle, we may interpret the
infinitesimal action of a meromorphic function $f$ deforming $\cF$
as changing the transition functions of $\cF$ by scalar
multiplication by $f$ on the locus where $f$ is defined. We may also
work with a formal variant, deforming sheaves using Laurent series
along a divisor in a smooth variety. For example, if $X$ is a curve,
$x\in X$ and $\K_x\supset \Oo_x$ are Laurent and Taylor series at
$x$ we have surjections
$$\K_x\twoheadrightarrow \K_x/\Oo_x=\HH^1_x(X,\Oo_X)\twoheadrightarrow
\HH^1(X,\Oo_X)$$ from $\K_x$ to local to global cohomology of $\Oo$,
which we can use to construct deformations of sheaves.

More generally, take a local section of the local cohomology sheaf
$\HH^1_D(Y,\Oo)$, where $D\subset Y$ is a divisor (locally principal
subscheme).
 This section
corresponds to an element of
$$\Hom(\Oo,\Oo(\infty D)/\Oo) = \Hom(\Oo,
\widehat{\Oo}(\infty D)/\widehat{\Oo})$$ where we pass to
completions along $D$. Tensoring this homomorphism by a sheaf $\cF$
we obtain a homomorphism $$\cF \to \widehat{\cF}(\infty
D)/\widehat{\cF}.$$ Now we pull back the canonical extension
$$0\to\cF \to \cF(\infty D) \to \widehat{\cF}(\infty D)/\widehat{\cF}\to 0$$
along this map and to obtain the desired local extension of $\cF$ by
itself.

Thus we can deform sheaves on a curve using principal parts of
functions at a point. We refer to this construction as ``tweaking"
of sheaves. The flows of many algebraically integrable systems can
be described in this fashion (see e.g. \cite{DM} where the
Heisenberg flows of the KP hierarchies are described in this way and
\cite{BF} for the case of generalized Drinfeld-Sokolov
hierarchies). In the next section we explain how the flows of the
$GL_n$ Hitchin system, which are generically given by the action of
Picard groups of spectral curves which live in the cotangent bundle
$T^*X$ of a curve, are in fact uniformly given by the action of the
global cohomology $\HH^1(T^*X,\Oo)$ of the cotangent bundle.

\subsubsection{Tweaking Algebroids}\label{tweaking algebroids}
In this section we consider a more general construction,
constructing arbitrary deformations of sheaves near a divisor
$D\subset Y$. This is modeled on the loop algebra uniformization of
moduli of bundles on a curve. The resulting flows do not commute in
general and form an action of a Lie algebroid on moduli spaces of
sheaves, rather than of a fixed Lie algebra (in other words, the
space parameterizing deformations depends on the sheaf being
deformed).

Let $\cE$ denote the $\Oo_Y$-algebra of Laurent series along the
divisor $D$, i.e. functions on the punctured formal neighborhood of
$D$. More precisely, $\cE$ is the inductive limit
$$\cE=\lim_{\longrightarrow}\; \wh{\Oo}_{Y,D}(kD),$$
where $\wh{\Oo}_{Y,D}$ is the completion of $\Oo_Y$ along $D$.

Consider a pair $(\cF,\xi)$ consisting of a coherent sheaf $\cF$ on
$Y$ and a Laurent endomorphism $\xi\in\End_{\cE}(\cF_{\cE})$ (where
$\cF_{\cE}=\cF\ot_{\Oo_Y}\cE)$) of the restriction of $\cF$ to the
punctured neighborhood of $D$. We then construct a first-order
deformation of $\cF$, $[\xi]\in \Ext^1(\cF,\cF)$, as the image
(under a connecting homomorphism) of the operation of restricting
sections of $\cF$ to $\cF_\cE$:
$$\{s\mapsto \xi(s_\cE) \mbox{ mod }\cF\}\in\Hom(\cF,\cF_{\cE}/\cF)\to \Ext^1(\cF,\cF).$$
More geometrically (and informally), the deformation is defined by
changing the transition function between the restrictions
$\cF|_{Y\sm D},\cF_{\wh{\Oo}_{Y,D}}$ of $\cF$ to $Y\sm D$ and to the
completion along $D$: we define a new sheaf by multiplying the
isomorphism between the restrictions of the two sheaves to the
punctured formal neighborhood (which are both $\cF_\cE$) by
$1+\epsilon \cdot \xi$ over the dual numbers
$\C[\epsilon]/\epsilon^2$.

In the case of vector bundles on a curve $X$ with $D=x\in X$, this
formal deformation procedure becomes the action of twisted loop
algebras at $x$ on the moduli space of vector bundles on $X$. If we
trivialize a vector bundle $\cF$ near $x$ then $\xi$ becomes an
element of the loop algebra $L\gl_n=\gl_n\ot\cE$ of Laurent series
of matrices at $x$, which acts on the moduli of bundles by
infinitesimally changing the transition functions at $x$. Without
the choice of trivialization, the twisted loop algebras
$\End_{\cE}(\cF_{\cE})$ form a transitive Lie algebroid over the
moduli stack of bundles.

Similarly in the general setting we may consider all tweakings
$\xi\in\End_{\cE}(\cF_{\cE})$ as forming a Lie algebroid over the
moduli stack $\Mf$ of coherent sheaves $\cF$ on $Y$. The tangent
sheaf to $\Mf$ is the sheaf $\ul{\Ext}^1$ of self--extensions along
$Y$ of the universal sheaf $\ul{\cF}$ on $\Mf$.

\begin{defn}
\mbox{}
\begin{enumerate}
\item The {\em algebroid of tweakings along $D\subset Y$} is the sheaf
$\ul{\End}_{\cE}=\End_{\cE}(\ul{\cF}_\cE)$ on $\Mf$ of $\cE$-module
endomorphisms of the universal sheaf. The anchor map is the map
$\ul{\End}_{\cE}\to\ul{\Ext}^1$ defined above.

\item The {\em algebroid of central tweakings} along $D$ is the image of
$\cE$ in $\ul{\End}_{\cE}$.
\end{enumerate}
\end{defn}

The central tweakings are the multiplication operators by functions
on the support of a sheaf $\cF_{\cE}$. Note that if we consider
sheaves where $\cF_\cE$ is a line bundle on its support, then all
tweakings are central: $\End_{\cE}(\cF_\cE)$ is given by functions
$\cE|_{\on{Supp}\cF_\cE}$ on the support of $\cF_\cE$. Thus, in this
generic case the Lie algebroid reduces to the commutative tweaking
action of meromorphic functions considered above. This is, in
particular, the case for CM spectral sheaves with simple framing,
i.e., $\framing=\bigoplus_1^k\Oo_{x_i}$ is a direct sum of
skyscrapers at distinct points of $E_\infty$ (for example in the
spinless case $\framing =\Oo_b$). For general spectral sheaves, we
obtain instead a richer nonabelian hierarchy of flows given locally
by the action of several copies of $L\gl_k$ for different $k$.

\subsubsection{Tweaking Framed Sheaves}\label{framed tweaking}
It is useful to consider also a relative version of the above
constructions for sheaves framed along a divisor $D$. Namely, we
would like to deform sheaves with a fixed restriction to $D$. First,
we have an action of the group-scheme of line bundles equipped with
a trivialization along $D$, or on the infinitesimal level of
$\HH^1(Y,\Oo(-D))$, lifting the action of $\HH^1(Y,\Oo)$ on
underlying unframed sheaves. More generally, meromorphic germs of
functions along $D$ act on framed sheaves in the same way as on
unframed sheaves (in the latter case the action depends on the germ
up to germs regular on $D$, while in the former germs are taken up
to those vanishing along $D$). More precisely, consider a coherent
sheaf $T$ on $D$ and the moduli stack $\Mf(T)$ of coherent sheaves
$\cF$ on $Y$ with an isomorphism $\cF|_D\to T$. We have a forgetful
map $\Mf(T)\to \Mf$ to the moduli stack of the underlying unframed
sheaves. We then have the following obvious lifting of the tweaking
algebroid:

\begin{lemma}
The pullback of the sheaf $\ul{\End}_{\cE}$ of tweakings along $D$
from $\Mf$ to $\Mf(T)$ has a canonical structure of Lie algebroid
lifting the action on $\Mf$. The anchor map on $\Mf$ vanishes on
endomorphisms regular on $D$, while that on $\Mf(T)$ vanishes on
endomorphisms vanishing on $D$.
\end{lemma}

\subsection{Flows of Hitchin Systems}\label{Hitchin flows
section} In this section we give an explicit description of the
correspondence between the Hitchin hamiltonians and their flows on
moduli spaces of Higgs bundles. This description may be viewed as a
geometric (or spectral) reformulation of the Lax pairs with spectral
parameter for hamiltonian flows on loop algebras (see
\cite{AHP,BBT,DM, LiMulasePrym, LiMulaseHitchin}). The technique is based on the loop group
uniformization of moduli of bundles, and parallels the discussion of
isomonodromy flows in \cite{sug}.

We start with a trivial statement about the prototype for the
construction, the basic hamiltonian system on $T^*GL_n$ (the
geodesic flow on $GL_n$). Let us identify $\gl_n^*$ with $\gl_n$ by
the trace form and $T^*GL_n$ with $GL_n\times \gl_n$ by left
translation. For every positive integer $i$, we have a function on
$T^*GL_n$ given by $\frac{1}{i+1}\on{Tr}Y^{i+1}$ on the $\gl_n$
component. There are two natural ways to write the corresponding
hamiltonian vector field. In Lax form, we describe the vector field
at $(X,Y)\in GL_n\times\gl_n$ as $([Y^i,\cdot],0)$. On the other
hand we have a spectral interpretation of this flow. We will
identify the adjoint quotient stack $\gl_n/GL_n$ with the space of
torsion coherent sheaves on the affine line $\spec \C[t]$ of length
$n$. The Lie algebra $\gl_n$ itself is identified with $\C[z]$
module structures on $\C^n$, i.e. torsion sheaves as above with a
basis. Finally $GL_n\times \gl_n$ is identified with torsion sheaves
with two bases, in other words module structures on $\C^n$ plus an
additional basis for $\C^n$. In this formulation, the hamiltonian
vector field above simply rescales the second basis by the action of
the function $z^i$ on $\spec\C[z]$.

We would like to have an analogous spectral description of Hitchin
hamiltonian flows, by reinterpreting the corresponding Lax form. Let
$X$ denote a smooth connected projective curve, and $Bun_n(X)$ the
moduli stack of rank $n$ vector bundles on $X$. Fix a point $x\in
X$, and let $\K\simeq \C((z))$ and $\Oo\simeq\C[[z]]$ denote the
complete local field and ring at $x$ (functions on the punctured
disc $D^\times$ and disc $D$ at $x$, respectively). We denote by
$Bun_n(X,x)$ the moduli scheme of bundles with a trivialization on
$D$, i.e. an infinite level structure at $x$. The cotangent bundles
$T^*Bun_n(X)$ and $T^*Bun_n(X,x)$ are identified with the moduli of
rank $n$ Higgs bundles $\eta\in \HH^0(X,\End V\ot \Omega)$ ($V\in
Bun_n(X)$) on $X$ and the moduli of Higgs bundles
$\eta\in\on{H}^0(X\sm x,\End E\ot \Omega)$ ($V\in Bun_n(X,x)$) with
arbitrary pole at $x$, respectively.

We will identify these cotangent bundles with moduli spaces of
spectral sheaves on $T^*X$ \cite{Hi,M,DM}---namely a Higgs field
$\eta$ on $V\in Bun_n(X)$ gives $V$ an $\Oo_{T^*X}$-module
structure, defining a coherent sheaf on $T^*X$ which pushes forward
to $V$. Likewise a meromorphic Higgs field $\eta$ as above makes
$V(\infty\cdot x)=j_*j^* V$ (where $j$ is the inclusion of $X\sm x$
into $X$) into an $\Oo_{T^*X}$-module. The support of this sheaf
(the {\em spectral curve}) is the zero locus of the characteristic
polynomial of $\eta$, i.e. the spectral curve is $X_\eta=\spec_X
\Oo_{T^*X}/\{\on{char}\eta\}$. We will describe the Hitchin flows on
this data in the same way as in the prototype example (with the
Higgs field as an $\Omega$--twisted version of the matrix $Y$ over
the punctured disc).

The spectral sheaf interpretation of the Hitchin system allows us to
define infinitesimal actions of $\HH^1(T^*X,\Oo)$ on $T^*Bun_n(X)$
and of meromorphic germs along the cotangent fiber $F_x\subset T^*X$
on $T^*Bun_n(X,x)$. The latter is defined by tweaking the
corresponding spectral sheaf as before, and preserving the
trivialization of $V|_{D}$ (note that this is well-defined since the
tweaking flows canonically preserve $V|_D$ while changing the gluing
with $V|_{X\sm x}$).

 Let
\begin{eqnarray*}
\Hitch(X)&=&\bigoplus_{n=1}^\infty\HH^0(X,\Omega^{i+1})\hskip.1in \subset\\
\Hitch(X,x)&=& \bigoplus_{i=0}\HH^0(X\sm x,\Omega^{\ot{i+1}})
\hskip.1in \subset\\ \Hitch(\Dx)&=&\bigoplus_{i=0} \Omega_{\K}^{i+1}
\end{eqnarray*} be the infinite Hitchin base spaces (where $\Omega_\K\simeq \C((z))dz$ is the
space of Laurent differentials at $x$). The Hitchin maps
$$H:T^*Bun_n(X)\to \Hitch(X),\hspace{1em} H: T^*Bun_n(X,x)\to
\Hitch(X,x)$$
have as the $i$th component $$H_i:(V,\eta)\mapsto
\frac{1}{i+1}\on{Tr}(\eta^{i+1}).$$ Note that we write the Hitchin
map in the basis for invariant polynomials coming from traces of
powers, rather than the usual basis consisting of coefficients of
the characteristic polynomials.

The topological dual of the (Tate) vector space $\K$ is identified
by the residue pairing with the differentials $\Omega_\K$. Likewise
if $\cT_\K\simeq \C((z))\del_z$ is the space of Laurent vector
fields, then $\cT_{\K}^i$ is identified with the dual to
$\Omega_\K^{i+1}$. Let $$\Oo(T^*\Dx)=\bigoplus \cT_\K^i$$ denote
functions on the (suitably defined) cotangent bundle of
$\Dx$.

\begin{lemma}
The graded dual spaces of the Hitchin spaces are canonically
identified as follows:
\begin{eqnarray*}
\Hitch(\Dx)^*&=& \Oo(T^*\Dx)\hskip.1in
\twoheadrightarrow \\
\Hitch(X,x)^*&=& \Oo(T^*\Dx)/\Oo(T^*(X\sm x)) \hskip.1in \twoheadrightarrow\\
 \Hitch(X)^*&=&\HH^1(T^*X,\Oo)
\end{eqnarray*}
\end{lemma}

It follows that we may identify $\HH^1$ classes on $T^*X$ with
linear functions on the Hitchin base space, and meromorphic germs on
$T^*X$ along the fiber $F_x$ with linear functions on the
meromorphic Hitchin space.

\begin{thm}\label{Hitchin post}
For a class $\xi\in \HH^1(T^*X,\Oo)$ (respectively
$\xi\in\Oo(T^*\Dx)$), the hamiltonian vector field on $T^*Bun_n(X)$
($T^*Bun_n(X,x)$) associated to $H^*\xi$ ($H^*_x\xi$) is identified
with the tweaking action of $\xi$.
\end{thm}

\begin{proof}
Morally, the theorem follows by hamiltonian reduction from the
corresponding statement on $T^*GL_n(\K)$. More concretely, let
$(V,\eta)\in T^*Bun_n(X,x)$. The tangent space to the moduli space
at $(V,\eta)$ is given by $$(s,\theta)\in \End_{\Dx}(V)/\End_{X\sm
x}(V) \oplus \HH^0(X\sm x,End(V)\ot\Omega).$$ The symplectic form on
the tangent space is given by the residue of the trace pairing:
$$\omega\left((s_1,\theta_1),(s_2,\theta_2)
\right) =\on{Res}(\on{Tr}(s_1\theta_2-s_2\theta_1)).$$ Fix $\xi\in
\cT_{\K}^i\subset\Oo(T^*\Dx)$ (we assume $\xi$ homogeneous for
simplicity of notation). The corresponding Hitchin hamiltonian is
$H_\xi(V,\eta)=\on{Res}(\xi\on{Tr}(\frac{1}{i+1}\eta^{i+1}))$.
Perturbing $(V,\eta)$ by $(s,\theta)$ as above, we find the
differential of this function is
$$dH_{\xi}|_{(V,\eta)}(s,\theta)=\on{Res}(\xi\on{Tr}(\eta^i\theta)).$$
To calculate the hamiltonian vector field
$v_{\xi}|_{(V,\eta)}=(s_{\xi},\theta_{\xi})$,  we must solve
$$\on{Res}(\xi\on{Tr}(\eta^i\theta))= dH_\xi|_{(V,\eta)}(s,\theta) =
\omega_{(V,\eta)}(v_\xi, (s,\theta)) =
\on{Res}(\on{Tr}(s_{\xi}\theta)) -\on{Res}(\on{Tr}(s\theta_\xi)).$$
It is immediate that
 $v_\xi|_{(V,\eta)}=(\xi\eta^i,0)$. In other words, we have
written the Hitchin flows in Lax form, with the flow at $(V,\eta)$
given by the action on $V$ of the element $\xi\eta^i\in \on{End}(V|_{D^\times})$.

On the other hand the tweaking action of a homogeneous
$\xi\in\cT_\K^i$ is given by multiplication by $\xi\in\Oo(T^*\Dx)$
as an endomorphism of $V|_{\Dx}$, considered as an
$\Oo(T^*\Dx)$-module using $\eta|_{\Dx}$. However this endomorphism
is simply the product $\xi\eta^i\in\End_{\Dx}(V)$, and the tweaking
action by this endomorphism is precisely the Lax vector field we
derived above. Thus the (central) tweaking flow by $\xi$ on $T^*X$
is written as an element of the tweaking algebroid at $x$ on $X$,
where it is identified with the loop algebra action on bundles $V$
with trivialization on the disc.

Finally note that if $\xi\in \Oo(T^*(X\sm x))\subset \Oo(T^*\Dx)$,
then on the one hand the Hitchin hamiltonian defined by $\xi$
vanishes on $T^*Bun_n(X)$, while on the other the tweaking action of
$\xi$ on spectral sheaves is trivial, since it is given by a global
rescaling on $X\sm x$ which vanishes in $\Ext^1$. This establishes
the theorem in the meromorphic case.

For $(V,\eta)\in T^*Bun_n(X)$, we choose a trivialization of $V|_D$
to put ourselves in the previous situation. The Hitchin polynomials
of $\eta$ are now regular on $D$, so $\xi\in \Oo(T^*\D)\subset
\Oo(T^*\Dx)$ pull back to the zero function of $(V,\eta)$. Likewise,
the tweaking action of such $\xi$ corresponds to a regular
endomorphism of $V$ on $T^*\D$, which corresponds to a change of
trivialization of $V$ on $D$. Thus the action on $T^*Bun_n(X)$
descends to
$$\Oo(T^*(X\sm x))\backslash \Oo(T^*\Dx)/ \Oo(T^*D)=\HH^1(T^*X,\Oo),$$
as desired.

\end{proof}

\begin{remark} For simplicity, we have stated the theorem in the setting
of a smooth curve $X$. However, the result carries over to a more
general setting. In particular we are interested in the case when
$X=E$ is a cubic curve, $x\in \G$ is a smooth point and $T^*E$ is
the log cotangent bundle. It is evident that the above arguments
(which are local at $x$) extend automatically to this
setting.
\end{remark}

\subsection{Compatibility of Hitchin and Heisenberg Flows}
In \cite{DM} and \cite{LiMulaseHitchin}, the authors investigate the
compatibility of the Hitchin systems on $T^*Bun_n(X,D)$  with KP- or
KdV-type flows on moduli of spectral curves defined using the
Krichever construction. More precisely, if we fix a point $x\in X$
and a partition $p$ of $n$, we can look at the moduli stack
$\on{Higgs}_{x,p}\to T^*Bun_n(X,D)$ of Higgs bundles whose spectral
curve over the punctured neighborhood of $x$ is identified with a
fixed spectral cover $\Sigma_p\to D^\times$ over the punctured disc,
with ramification type $p$. The partitions $p$ label conjugacy
classes of Heisenberg (or Cartan) subalgebras of the loop algebra
$L\gl_n$,consisting of loops into diagonalizable matrices whose
eigenvalues undergo a permutation of type $p$ around the puncture.
By fixing a distinguished ramified cover $\Sigma_p$ of the punctured
disc of type $p$, we pick out a particular Heisenberg algebra
$A_p\subset L\gl_n$, isomorphic to functions on $\Sigma_p$, and
there is a natural action of $H_p$ on $\on{Higgs}_{x,p}$ by
tweakings, modifying the Higgs bundle but preserving its spectral
curve (see \cite{AB} and \cite{BF} for detailed discussions of
Heisenberg and KdV flows). The above papers prove that the actions
of these Heisenberg algebras are hamiltonian with respect to the
Poisson structure on Higgs bundles, when sufficiently strong
conditions are imposed on the regularity of the spectral curves or
the Higgs fields at $x$. Since the $A_p$ action is given by
tweakings, we immediately recover from Theorem \ref{Hitchin post} a
strong form of the compatibility, independent on the singularities
of the spectral curve. (Note that on $\on{Higgs}_{x,p}$ all spectral
sheaves are rank one on their support near $x$, so all tweakings
along $F_x$ are central.)

\begin{corollary}[Compatibility]
\label{compatibility}
Fix $x\in X$ and a partition $p$ of $n$. The algebroid of tweakings
along $F_x\subset T^*X$ pulls back to the action of the Heisenberg
algebra $A_p$ on $\on{Higgs}_{x,p}$. Therefore all Heisenberg flows
are given by pullbacks of Hitchin hamiltonians.
\end{corollary}

\subsection{Framed Calogero-Moser Hierarchies}\label{CM flows}

In this section we consider the tweaking flows on CM spectral
sheaves. Since we are considering sheaves on $\Enatbar$ framed along
$E_\infty$, a natural collection of tweaking flows is parametrized
by the algebroid of tweakings along $E_\infty$:

\begin{defn}
The {\em CM algebroid} on $\CM_n(E,\framing)$ is the Lie algebroid
$\ul{\End}_\cE$ of tweakings along the divisor $E_\infty\subset
\Enatbar$.
\end{defn}

We would like to obtain concrete hierarchies of tweakings,
specifically tweaking by meromorphic functions (central tweakings).
To do so we single out particular meromorphic germs on $\Enatbar$ by
which to tweak. We will then identify the resulting flows with
Hitchin hamiltonian flows on Higgs bundles.

 Let us consider first the $k$-spin case
$\framing=\Oo_b^k$. Recall the meromorphic function $\underline{t}$
on $\Enat$ defined near the fiber $F_b$ over $b$ (Section \ref{about
Enat}). We will define flows on CM spectral sheaves as tweaking near
$b\in E_\infty$ by powers of $\ul{t}$. Note that the definition of
$\underline{t}$ required a choice of global vector field on $E$,
equivalently a trivialization of the cotangent fiber of $E$.
Accounting for the action of change of trivialization, we identify
the resulting polynomial algebra $\C[\ul{t}]$ canonically with
$\C[\del]$ where $\del$ denotes a global vector field on $E$.

Given a polynomial $p$ in $\C[\ul{t}]$ (considered as a germ of a
meromorphic function near the fiber $F_b$) and a CM spectral sheaf,
we may restrict $p$ to define a germ of a meromorphic function on
the corresponding CM spectral curve at its intersection point
$b_\infty$ with $E_\infty$. This defines an $\HH^1$ class of
$\Oo(-E_\infty)$, which then acts as a deformation of the framed
spectral sheaf. Thus we have an action of $\C[\ul{t}]$ by commuting
vector fields on $\CM_n^k(E)$:

\begin{defn}[CM flows]\label{CM tweaking}
The {\em spin CM hierarchy} on spectral sheaves is the canonical action of
$\C[\del]$ on $\CM_n^k(E)$, where $\del^i$ acts by deforming
spectral sheaves by the restriction of $\ul{t}^i$.
\end{defn}

We may similarly define CM flows on $\framing$-framed spectral
sheaves for an arbitrary $\framing$. Let $S=\coprod_{j=1}^l
x_j\subset E$ denote the set-theoretic support of $\framing$,
consisting of $l$ distinct points on $X$ (counted without
multiplicity).

Let $\wh{F}_j$ denote the formal neighborhood of the fiber
$F_{x_j}\subset \Enat$. For $x_j\neq b$, we consider the restriction
of the polynomial ring $\C[\del]$ of global functions on $T^*E$ to
$\wh{F}_j$, giving meromorphic endomorphisms of framed spectral
sheaves. When $x_j=b$, we shift these functions by the zeta function
as above, again giving an action of $\C[\del]$ (as powers of
$\ul{t}$) as meromorphic endomorphisms of framed spectral sheaves
restricted to $\wh{F}_b$.

\begin{defn}\label{framed CM hierarchy def}
The {\em framed CM hierarchy} on $\CM_n(E,\framing)$ is the action of
$\C[\del]^l$ on $\CM_n(E,\framing)$, where $\del^i$ in the $j$th
copy acts by tweaking on $\wh{F}_j$ by $\del^i$ for $x_j\neq b$ and
by $\ul{t}^i$ for $x_j=b$.
\end{defn}

\begin{thm}\label{framed CM flows thm}
\mbox{}
\begin{enumerate}
\item The identification $\CM_n(E,\framing)\to\HCM_n(E,\framing)$ of
Theorem \ref{CM identification} (in particular $\CM_n^k(E) \to
\HCM_n^k(E)$) intertwines the framed CM hierarchy (action of
$\C[\del]^l$) with the flows of the framed CM hamiltonians,
Definition \ref{framed CM hamiltonians}: $\del^i$ at the point $x_j$
is the hamiltonian flow for $H_{i,x_j}$.

\item For simple framing $\framing =\bigoplus_{i=1}^k\Oo_{x_i}$ this construction gives all
tweaking flows: the action of the full CM algebroid is identified
with the action of $\C[\del]^k$.
\end{enumerate}
\end{thm}

\begin{proof}
By Theorem \ref{Hitchin post}, the Hitchin flows are defined by
tweaking of spectral sheaves on $T^*E$, while the CM spectral flows
are defined by tweaking sheaves on $\Enat$ framed along $E_\infty$.
More specifically, the hamiltonian flow of the spin CM hamiltonian
$H_{i,b}$, given by a residue of a trace of $(\eta+\zeta\on{Id})^i$,
is identified by Theorem \ref{Hitchin post} with the tweaking action
of the function $(\del+\zeta)^i$ on $T^*E$, which has a pole along
the fiber $F_b$ over $b$. This function is identified with
$\ul{t}^i$ by the birational identification of $T^*E$ and $\Enat$.
It follows that under the isomorphism of Corollary \ref{CM
identification}, the meromorphic endomorphism of a Higgs bundle
given by $(\del+\zeta)^i$ is identified with the meromorphic
endomorphism of a twisted Higgs bundle, hence CM spectral sheaf,
given by $\ul{t}^i$. For $x_j\neq b$, our hamiltonians are simply
residues of traces of powers of the Higgs field, which correspond to
multiplication by powers by the coordinate function $\del$ on the
fibers of $T^*E$ (which is identified with $\Enat$ away from $F_b$).
(Note that the tweaking flows preserve the framing at $E_\infty$,
while the Hitchin hamiltonians $H_i$ are independent of the framing
data $u,v$.) It follows that the corresponding tweaking flows are
identified, as claimed.

Part (2) follows from the observation (Section \ref{tweaking
algebroids}) that for simple framing, the $T$-framed spectral
sheaves $\cF$ are all rank one torsion-free sheaves on their
support. Hence all tweakings along $E_\infty$ are given by
multiplication by meromorphic functions on the support of the
microlocalization $\cF_\cE$, and considered up to functions that
vanish on the intersection with $E_\infty$. Finally the support of
$\cF_\cE$ is a union of $k$ punctured discs $\Dx_j$ (one through
each $x_j\in E_\infty$). Hence functions on $\Dx_j$ are identified
(via restriction) with Laurent series $\C((\del\inv))$ in the fiber
coordinate $\del$ on $T^*E$ (or $\ul{t}$ for $x_j=b$). The assertion
follows from the observation the polynomial algebra $\C[\del]$ maps
isomorphically to the quotient of $\C((\del\inv))$ by functions
vanishing at $x_j$.

\end{proof}

\bibliographystyle{alpha}
\newcommand{\bl}{/afs/math.lsa.umich.edu/group/fac/nevins/private/math/bibtex/}

\bibliography{\bl baranovsky,\bl bradlow,\bl brosius,\bl
bialynicki-birula, \bl carrell,\bl donagi,\bl friedman,\bl
frankel,\bl ginzburg,\bl gomez,\bl griffiths, \bl grojnowski,\bl
hartshorne,\bl hitchin,\bl hungerford,\bl huybrechts,\bl jardim,\bl
kapranov,\bl kapustin,\bl kirwan, \bl kurke,\bl langton,\bl
laumon,\bl lehn,\bl matsumura,\bl milne, \bl nakajima,\bl nevins,\bl
segal,\bl sernesi,\bl thaddeus,\bl viehweg,\bl wilson}

\end{document}